\documentclass[12pt]{amsart}

\usepackage{amsmath}
\usepackage{latexsym}
\usepackage{amssymb}
\usepackage{amsfonts}
\usepackage{xcolor}

\usepackage{graphics}

\include{PDF}

\renewcommand{\proof}{\par\noindent{\it Proof.\ \ }}
\def\qed{\ifmmode\square\else\nolinebreak\hfill
$\Box$\fi\par\vskip12pt}

\def\ov{\overline} 
\def\l{\langle} \def\r{\rangle}

 \def\ZZ{{\mathbb Z}}

\def\calM{{\mathcal M}}

 \def\MM{{\mathcal M}} 
\def\ZZ{{\mathbb Z}}

\def\K{{\bf K}}

\def\Aut{{\rm Aut}}

\def\Cos{{\rm Cos}}
\def\D{{\rm D}} 
\def\S{{\rm S}}

\def\C{{\bf C}}\def\Z{{\bf Z}}

\def\Ga{{\it \Gamma}}
\def\Ome{{\it \Omega}}
\def\Sig{{\it \Sigma}}
\def\Del{{\it \Delta}}

\def\a{\alpha} \def\b{\beta} \def\g{\gamma} 
 
\def\o{\omega}

\def\A{{\rm A}}\def\Sym{{\rm Sym}}

  \def\D{{\rm D}}

\newtheorem{theorem}{Theorem}[section]%
\newtheorem{lemma}[theorem]{Lemma}%
\newtheorem{corollary}[theorem]{Corollary}%
\newtheorem{proposition}[theorem]{Proposition}%
\newtheorem{definition}[theorem]{Definition}%
\newtheorem{example}[theorem]{Example}%
\newtheorem{construction}[theorem]{Construction}%
\newtheorem{remark}[theorem]{Remark}%
\newtheorem{hypothesis}[theorem]{Hypothesis}%

\renewcommand{\proof}{\par\noindent{\it Proof.\ \ }}
\def\qed{\ifmmode\square\else\nolinebreak\hfill
$\Box$\fi\par\vskip12pt}

\def\l{\langle} \def\r{\rangle}
 
\renewcommand{\leq}{\leqslant}
\renewcommand{\geq}{\geqslant}
 
\def\ov{\overline} 

 \def\bfH{{\bf H}} \def\bfI{{\bf I}}

 \def\bfQ{{\bf Q}}

\def\calM{{\mathcal M}}  
  
\def\calS{{\mathcal S}}

\def\MM{{\mathcal M}}

 \def\ZZ{{\mathbb Z}}

  \def\Cos{{\sf Cos}}
\def\Aut{{\sf Aut}}

\def\lcm{{\sf lcm}}

\def\Ga{{\it \Gamma}} \def\Ome{{\it \Omega}}\def\Sig{{\it \Sigma}}
\def\Del{{\it\Delta}}

\def\a{\alpha} \def\b{\beta} \def\o{\omega} \def\g{\gamma}

\def\C{{\bf C}}  

\def\A{{\rm A}} \def\S{{\rm S}}
 \def\D{{\rm D}}

\def\A{{\rm A}}
\def\Sym{{\rm Sym}} \def\Alt{{\rm Alt}}

   \def\D{{\rm D}}

\def\K{{\bf K}} 
 
\def\Z{{\bf Z}}

\def\RegMap{{\sf RegMap}}

\def\BRM{{\sf BiRoMap}}
\def\FRM{{\sf RegMap}}
\def\RM{{\sf RotaMap}}
\def\Core{{\sf Core}}
\def\SimpCos{{\sf SimpCos}}

\begin{document}

\title[Rotary maps]{Locally Finite Vertex-Rotary Maps and Coset Graphs with Finite Valency and Finite Edge Multiplicity}
\thanks{This work was partially supported  by NSFC grants 11701497, 11771200, 11931005
 and 61771019, and Australian Research Council Discovery Project grant DP160102323.}

\author{Cai Heng Li}
\address{SUSTech International Center for Mathematics, and Department of Mathematics\\
South University of Science and Technology of China\\
Shenzhen, Guangdong 518055\\
P. R. China}
\email{lich@sustech.edu.cn}

\author{Cheryl E. Praeger}
\address{Department of Mathematics and Statistics\\
The  University of Western Australia\\
Australia}
\email{cheryl.praeger@uwa.edu.au}

\author{Shujiao Song}
\address{Department of Mathematics and Statistics\\
Yantai University\\
P. R. China}
\email{shujiao.song@ytu.edu.cn}


\date\today

\begin{abstract}
It is well-known that a simple $G$-arc-transitive graph can be represented as a coset graph for the group $G$.
This representation is extended to a construction of $G$-arc-transitive coset graphs
$\Cos(G,H,J)$  with finite valency and finite edge-multiplicity, where $H, J$ are stabilisers in $G$ of a vertex and incident edge, respectively.
Given a group $G=\l a,z\r$ with $|z|=2$ and $|a|$ finite, the coset graph $\Cos(G,\l a\r,\l z\r)$ is shown, under suitable finiteness assumptions,  to have exactly two different arc-transitive embeddings as a $G$-arc-transitive map
$(V,E,F)$, namely, a {\it $G$-rotary} map if $|az|$ is finite, and a {\it $G$-bi-rotary} map if $|zz^a|$ is finite.
The $G$-rotary map can be represented as a coset geometry for $G$, extending the notion of a coset graph. However the $G$-bi-rotary map does not have such a representation, and the face boundary cycles must be specified in addition to incidences between faces and edges. We also give a coset geometry construction of a flag-regular map $(V,E,F)$.
In all of these constructions we prove that the face boundary cycles are regular cycles which are simple cycles precisely when the given group acts faithfully on $V\cup F$.

\medskip\noindent
{\it Keywords:}  arc-transitive  graphs,  coset graphs  with  multiple  edges,  graph  embeddings,  arc-transitive maps, rotary maps.

\medskip\noindent
{\it Math. Subj. Class.:}  20B25, 05C25, 05C35

\end{abstract}

\maketitle

\section{Introduction}\label{s:intro}

\medskip

The graphs we study in this paper are, formally, triples $\Gamma = (V, E, \bfI)$, where
$V, E$  are sets, the elements of which are called \emph
{vertices} and \emph{edges}, respectively, and $\bfI\subseteq V\times E$ is a binary relation called
{\it incidence} such that, for each edge $e\in E$, there are exactly
two distinct vertices incident with $e$, say $\alpha\, \bfI\, e$ and $\beta\, \bfI\, e$.
Thus a graph $\Gamma$ has no `loops' but a pair of `adjacent' vertices $\alpha,\beta$ (that is, $\alpha\, \bfI\, e$ and $\beta\, \bfI\, e$ for some edge $e$) may be incident with more than one edge. If $E\ne \emptyset$ then we say that $\Gamma$ is \emph{non-empty}.
If for each vertex $\alpha$, the set $E(\alpha) :=\{ e\in E \mid \alpha\,\bfI\,e\}$ is finite, then we say that $\Gamma$ is \emph{locally finite}. Moreover if there is a constant $\lambda$ such that, for each pair of adjacent vertices, there are exactly $\lambda$ edges incident with both, then
we say that $\Gamma$ has \emph{edge-multiplicity} $\lambda$.
We call $\Gamma$ \emph{simple} if it has edge-multiplicity $1$. \emph{All graphs in this paper will be locally finite.}

For graphs $\Ga=(V,E,\bfI)$ and $\Gamma'=(V', E', \bfI')$, an \emph{isomorphism} from $\Gamma$ to $\Gamma'$ is a bijection
$g: V\cup E \rightarrow V'\cup E'$ which maps $V$ to $V'$, and $E$ to $E'$, and induces a bijection
$\bfI\rightarrow \bfI'$.
We say that $\Gamma$ and $\Gamma'$ are isomorphic if such an isomorphism exists.
In particular, if $\Gamma = \Gamma'$ then an isomorphism is called an \emph{automorphism}:
it  is a permutation  $g$ of $V\cup E$ which fixes $V$ and $E$ setwise and satisfies $\alpha\ \bfI\ e$ if and only if
$\alpha^g\, \bfI\, e^g$ for all $\alpha\in V, e\in E$.
The set $\Aut\Ga$ of automorphisms forms a group under composition, called the \emph{automorphism group} of $\Ga$.
For a subgroup $G\leqslant\Aut\Ga$, the graph $\Ga$ is called {\it $G$-vertex transitive}
if $G$ is transitive on the vertex set $V$.
Similarly $\Ga$ is {\it $G$-edge-transitive} if $G$ is transitive on the edge set $E$.
Each edge $e\in E$ corresponds to two \emph{arcs} of $\Gamma$, denoted $(\alpha,e,\beta)$ and $(\beta,e,\alpha)$, where $\alpha, \beta$ are the two vertices incident with $e$, and  $\Ga$ is
{\it $G$-arc-transitive} if the induced $G$-action  on the set of arcs is transitive.
Since an arc  $(\alpha,e,\beta)$ is completely determined by the edge $e$,
together with  either the
initial vertex $\alpha$ or the final vertex $\beta$, sometimes the arc is  abbreviated to $(\alpha,e)$ or $(e,\beta)$. Also if we wish to emphasise in our notation the vertices $\alpha,\beta$ incident with an edge $e$ we write the edge as $[\alpha,e, \beta]$ or $[\beta,e,\alpha]$.

For a graph $\Ga=(V,E,\bfI)$ and a vertex $\alpha\in V$, the set
$\Gamma(\alpha)$ of vertices adjacent to $\alpha$ is called the \emph{neighbourhood} of $\a$,
and its cardinality $|\Gamma(\alpha)|$  is called the \emph{valency} of $\a$. If $\Ga$ is locally finite and
both $|E(\alpha)|$ and $k:=|\Gamma(\alpha)|$ are independent of $\alpha$, then $\Ga$ is said to be \emph{regular}
of valency $k$. In particular such a graph $\Ga$ has constant edge-multiplicity, say $\lambda$, and $|E(\alpha)|=
k\lambda$. Note that  each locally finite, arc-transitive graph is regular.

\subsection{Coset graphs}\label{ss:cosetgraphs}
We introduce a new coset graph construction for these not necessarily simple graphs (see Construction~\ref{def-m-graph}), which characterises all locally finite arc-transitive graphs.  See  Theorems~\ref{thm:cosetgraph} and~\ref{thm:cosetgraph2} for more details. For a subgroup $X$ of a group $G$, the \emph{core} of $X$ in $G$ is $\Core_G(X)=\cap_{g\in G}X^g$ (the largest normal subgroup of $G$ contained in $X$), and if the core is trivial then we say that $X$ is \emph{core-free} in $G$.

\begin{theorem}\label{new-cos-rep}
Let $G$ be a group, with subgroups $H, J$ such that $H\ne G$, $|J:H\cap J|=2$, $|H:H\cap J|$ is finite, and $H\cap J$ is core-free in $G$. Then  $\Cos(G,G_\a,J)$ (as defined in Construction~$\ref{def-m-graph}$) is a non-empty, locally finite $G$-arc-transitive graph. Conversely each non-empty, locally finite $G$-arc-transitive  graph is isomorphic to $\Cos(G, H, J)$, where  $H, J$ are stabilisers in $G$ of an incident vertex and edge, respectively.
\end{theorem}

Construction~\ref{def-m-graph} generalises a construction   introduced by Sabidussi \cite{Sab} in 1964 for $G$-arc-transitive simple graphs: if  $\Ga=(V,E, \bfI)$ is a locally finite, $G$-arc transitive  simple graph, then  $\Ga$ can be represented as the {\it simple coset graph}  $\SimpCos(G,G_\a,G_\a gG_\a)$, for a chosen vertex
$\a\in V$ and group element $g\in G$ such that $\a^g$ is adjacent to $\a$ in $\Ga$, namely $V = [G:G_\a] = \{ G_\a x\mid x\in G\}$, and  $E = \{ \{G_\a x, G_\a y\} \mid yx^{-1}\in G_\a g G_\a \}$,  with incidence given by inclusion. We tie these two constructions together as follows.

Each locally finite graph $\Gamma=(V, E, \bfI)$ with constant edge-multiplicity $\lambda$
(in particular each edge-transitive or arc-transitive graph)  corresponds to a unique simple \emph{base graph} $\Sigma=(V, E', \bfI')$, where $E'$ is the set of unordered pairs of adjacent vertices from $V$, and $\bfI'$ is inclusion.
On the other hand, given a positive integer $\mu$ and a graph $\Ga$ with constant edge-multiplicity $\lambda$, we may construct a graph $\Ga^{(\mu)}$
with constant edge-multiplicity $\lambda\mu$, called the  \emph{$\mu$-extender} of $\Ga$,  by
replacing each edge  $e$ of $\Ga$ with $\mu$ edges, each incident to the two vertices
incident with $e$. If both $\Ga$ and $\Ga^{(\mu)}$ are $G$-arc-transitive, for some group $G$, then
$\Ga^{(\mu)}$ is called a {\it $(G,\mu)$-extender} of $\Ga$. In particular,
if $\Ga$ is $G$-arc-transitive, then $G$ induces an arc-transitive group of automorphisms of its base graph
$\Sigma$, and $\Ga$  is a $(G,\lambda)$-extender $\Sigma^{(\lambda)}$ of  $\Sigma$.

\begin{theorem}\label{cos-rep}
Let $G, K, J$, and $\Ga=\Cos(G, H, J)$ be a locally finite $G$-arc-transitive graph $(V, E, \bfI)$
as in Theorem~$\ref{new-cos-rep}$, with $H\cap J$ core-free, and let $g\in J\setminus (H\cap J)$
and $L=(H\cap H^g)\l g\r$. Then $J\leq L$, and  $\Ga$ has edge-multiplicity $|L:J|$ and base graph $\Cos(G, H, L)\cong \SimpCos(G, H, HgH)$. In particular, $\Ga$ is simple if and only if $J=(H\cap H^g)\l g\r$, and in this case $\Ga \cong \SimpCos(G, H, HgH)$.
\end{theorem}

\subsection{Graph embeddings}\label{ss:graphembeddings}

We now apply the coset graph theory from Subsection~\ref{ss:cosetgraphs}  to develop a theory of arc-transitive embeddings for connected,
locally finite graphs into surfaces (for our concept of a graph). \emph{Thus from now on we restrict our attention to connected locally finite graphs, that is to say, the valencies, edge-multiplicities, and face-lengths are all finite.}

The book \cite{Tucker-book}
is a good reference for the concepts we now introduce (especially Chapter 1.4; note that the graphs considered there are more general than those
we study as loops are permitted).
By a {\it surface} we mean a connected compact 2-manifold without boundary.
Any sufficiently small open neighbourhood of each point is an open disc, so each point of such a surface is an interior point of an open disc.
Surfaces can be defined as complex algebraic curves, and, as we shall see,
they can also be defined by graph embeddings.
Given a surface $\calS$ and a connected locally finite  graph $\Ga=(V,E,\bfI)$, a {\it $2$-cell embedding} of $\Ga$ in $\calS$
is a `drawing' of $\Ga$ on $\calS$ such that $\calS\setminus(V\cup E)$ is a disjoint union of open discs,
that is, of 2-{\it cells}. Each $2$-cell in $\calS\setminus(V\cup E)$ is called a {\it face}, and we denote the set of faces by $F$.
We sometimes refer to a $2$-cell embedding, defined in this way, as simply an \emph{embedding} of $\Ga$,
and we more often refer to it as a {\it map}, and denote it by $\calM = \calM(V, E, F)$.

For each face $f\in F$, let $E_F(f)$ denote the set of edges on the boundary of $f$.
We  require not only that $\Gamma$ is locally finite (so all
valencies and edge-multiplicities  are finite), but also, for each face $f$,  we require $E_F(f)$ to be finite;
$|E_F(f)|$ is called the \emph{face length} of $f$.
The map $\calM$ is said to be \emph{locally finite} when both the underlying graph is locally finite and all face lengths are finite.

To describe $\calM$ further we need to formalise the concept of a cycle in $\Gamma$. For  a subset $S\subset E$, we denote by $[S]$ the {\it edge-induced subgraph} on $S$, namely,
the edge set is $S$ and  the vertex set consists of all vertices incident with some edge in $S$.

\begin{definition}\label{def:cycles}
{\rm
For an integer $\ell\geq3$, a {\it cycle of length $\ell$} in a graph is a sequence $C=(e_1,e_2,\ldots,e_\ell)$ of $\ell$ pairwise distinct edges
of the form $[\a_{i-1},e_i,\a_i]$ for each $i$  (reading subscripts modulo $\ell$),
so in particular $\a_0=\a_\ell$, and we often write $e_\ell=e_0$.
It is possible that some of the $\alpha_i$  may be repeated.

If the $\alpha_i$ are pairwise distinct then the edge-induced subgraph $[C]$ (with vertex set $\{\alpha_i\mid 1\leq i\leq \ell\}$ and edge set $\{e_i\mid 1\leq i\leq \ell\}$) is a simple graph, called a \emph{simple cycle} and denoted
$[C]=\C_\ell$.
If $\ell=n\lambda$, and if $[C]$ has $n$ distinct vertices such that $\a_i=\a_{i+n}$ for all $i$, reading the subscripts modulo $n$, then $n\geq2$: if $n\geq3$ then $[C]$ is the $\lambda$-extender $\C_n^{(\lambda)}$, while if $n=2$ then $\ell=2\lambda\geq 4$ and $[C]$ is the $\ell$-extender $\K_2^{(\ell)}$. In either of these two cases,
$C$ is said to be a \emph{regular cycle}.
}
\end{definition}

Since each face $f$ of a map  $\calM = \calM(V, E, F)$ is an open disc, the edge set $E_F(f)$  admits a natural sequencing
as a cycle of $\Gamma$ and this cycle may start with any of the edges in $E_F(f)$ and may move in the `natural or reverse direction' around the disc. The set of all
such cycles forms a single sequence class $C(f)$ of cycles,
and $C(f)$, or any of the cycles in it, is  called the \emph{boundary cycle} of $f$.
We emphasise that \emph{a construction of a map must specify the sequence class $C(f)$ of boundary cycles for each face $f$, not just the edge-subset $E_F(f)$}.
Combinatorially it is sometimes useful to identify a face with its boundary cycle, that is, we sometimes simply write $f$ for the sequence class $C(f)$.

\begin{definition}\label{def:circembedding}
{\rm
Let $\calM = (V, E, F, \bfI)$ be a map with underlying graph $\Ga=(V, E, \bfI)$.
Then $\calM$ is called a {\it circular embedding} of $\Ga$ if $C(f)$ is a simple cycle  of $\Ga$,  for each face $f\in F$.
We also say that $\calM$ is a \emph{circular map}.
}
\end{definition}

The vertices and edges of $\Gamma$ on the boundary cycle of a face $f$ are defined to be \emph{incident} with $f$.
In this way the incidence relation $\bfI$ defining $\Gamma$ is extended to a binary relation on $V\cup E\cup F$,
including these incidences between vertices and faces, and between edges and faces, and
\emph{we view the map $\calM$ as an incidence configuration $(V, E, F, \bfI)$, and recall that $F$ is equipped with the set of specified boundary cycles
of the faces.}
We sometimes omit reference to the incidence relation $\bfI$ and write $\calM=(V,E,F)$.
For an embedding $\calM$ of a connected graph $\Ga=(V,E, \bfI)$ in a surface,
the graph $\Ga$ is called the {\it underlying graph} of $\calM$, and
the surface $\calS$ is called the {\it supporting surface}  of $\calM$.
When working with maps, an {\it arc} $(\alpha,e,\beta)$ of the underlying graph is often
denoted $(\alpha, e)$ (as discussed above) and so an arc of the map may be considered as
an incident pair consisting of a vertex and an edge.
A {\it flag}  is a triple consisting of a vertex, an edge, and a face, each pair of which is incident.

An {\it isomorphism} from a map $\calM_1=(V_1,E_1,F_1,\bfI_1)$ to a map $\calM_2=(V_2,E_2,F_2,\bfI_2)$ is a bijection
$g:\ V_1\cup E_1\cup F_1\to V_2\cup E_2\cup F_2$ which maps $V_1$ to $V_2$, $E_1$ to $E_2$, $F_1$ to $F_2$
(and boundary cycles of $\MM_1$ to boundary cycles of $\MM_2$),
and induces a bijection $\bfI_1\rightarrow \bfI_2$.
We write $\calM_1\cong\calM_2$.
In particular the restriction of $g$ to $V_1\cup E_1$ is an isomorphism of the  underlying graphs $\Ga_1 \rightarrow \Ga_2$,
where $\Ga_1=(V_1,E_1, \bfI_1)$ and $\Ga_2=(V_2,E_2, \bfI_2)$, such that for $f_1\in F_1$, $g$ maps the boundary cycle $C(f_1)$
(a sequence class of cycles of $\Gamma_1$) to the boundary cycle $C(f_1^g)$ of $f_1^g\in F_2$
(a sequence class of cycles of $\Gamma_2$).
If $\calM:=\calM_1=\calM_2$ (so $\Ga:=\Ga_1=\Ga_2$), then $g$ is called an {\it automorphism} of $\calM$.
The set $\Aut\MM$ of automorphisms of $\MM$ forms a group called the automorphism group of $\calM$.
It is a subgroup of $\Aut\Ga$, and the inclusion can be proper. If an automorphism $g$ of $\calM$ fixes a flag $(\a,e,f)$
of $\calM$, then $g$ fixes each vertex and edge incident with $f$, and $g$ fixes each face incident with $\alpha$,
and so on, so that, since the underlying graph $\Ga$ of $\MM$ is connected, $g$ fixes
every element of $V\cup E\cup F$, and hence $g$ is the identity.

\begin{definition}\label{def:rotary}
{\rm
(a) For a connected locally finite graph $\Ga=(V, E, \bfI)$, and $G\leqslant\Aut\Ga$,
we say that
$\Ga$ is {\it $G$-vertex-rotary} if $G$ is arc-transitive and the stabiliser $G_\alpha$ of a vertex $\alpha$ induces a transitive cyclic group on $E(\alpha)$; sometimes we simply say that $\Ga$ is {\it vertex-rotary} if such a group $G$ exists.

(b) For a locally finite map $\calM=(V, E, F, \bfI)$, and $G\leqslant\Aut\calM$, we say that $\MM$ is {\it $G$-arc-transitive} if $G$ is transitive on the arc set,
 that $\MM$ is {\it $G$-flag-regular} if $G$ is transitive (and hence regular) on the set of flags of $\calM$, and that
$\calM$ is {\it $G$-vertex-rotary} if the underlying graph $\Ga=(V, E, \bfI)$ is  $G$-vertex-rotary.
Similarly a $G$-arc-transitive map $\calM$ is called  {\it $G$-face-rotary} if the stabiliser $G_f$ of a face $f$ induces a
cyclic group on $E_F(f)$, and hence on the boundary cycle $C(f)$. The map $\calM$ is said to be {\it $G$-rotary} if it is both $G$-vertex-rotary and $G$-face-rotary.

(c) It is possible that a $G$-vertex-rotary map is not $G$-face-rotary; it turns out that for such maps, the local orientations of the carrier surface induced by the cyclic vertex stabilisers disagree as one traverses an edge, and these maps are called \emph{$G$-bi-rotary}, see the discussion in Remark~\ref{rem:birot}  and also see \cite{bi-rotary} where this terminology is introduced when the full automorphism group of the map is arc-regular.
}
\end{definition}

Possibilities for the induced subgraphs $[C(f)]$ of the boundary cycles for faces $f$ of $G$-arc-transitive maps are studied in \cite[Theorem 1.2]{Cycles}.
To help in our study of vertex-rotary maps, we show in Lemma~\ref{lem:rotpair} that a connected locally finite vertex-rotary graph is determined by a group with two generators, one of which is an involution.

\begin{definition}\label{def-rotarypair}{\rm
For elements $a, z$ of a group $G$, the ordered pair $(a,z)$ is called a \emph{rotary pair} if $G=\langle a,z\rangle$, $|a|$ is finite, and $|z|=2$ (that is, $z$  is an involution) with $z\not\in\l a\r$.
}
\end{definition}

Lemma~\ref{lem:rotpair} shows further that the connected $G$-vertex-rotary graph $\Ga=(V, E, \bfI)$ is isomorphic to the coset graph $\Cos(G,\l a\r,\l z\r)$ as defined in Construction~\ref{def-m-graph}, and that $G$ acts regularly on arcs. Moreover, $V$ can be identified with the set $[G:\l a\r]$ of right cosets of $\l a\r$, and $E$ with $[G:\l z\r]$ such that $\a=\l a\r$ is incident with $e=\l z\r$, so  $e=[\a,e,\b]$ where $\b=\a^z=\l a\r z$.
Thus $z$ fixes $e$ and interchanges the two vertices $\a$ and $\b$, and we have $G_\a=\l a\r$ and $G_e=\l z\r$.  Moreover, if $|a|=2$ then by Corollary~\ref{cor:rotpair}(c), either $\Gamma$ has just two vertices, or $\Gamma$ is a simple cycle and hence contains only one cycle; so any vertex-rotary map in these cases would have either at most two vertices, or at most two faces.
There are some interesting maps
$\calM=(V, E, F, \bfI)$ where either $|V|\leq2$ or $|F|\leq 2$.
Such maps are discussed in \cite{Li-Siran} for regular maps, and in \cite{Cycles2} for the arc-transitive case; and certain edge-transitive maps with a single face are described in \cite{Fan-Li,Singerman}.

In the light of these comments,  we  assume in our study of $G$-vertex-rotary maps that there are at least three vertices and at least three faces. This implies in particular that $|a|\geq 3$ for a rotary pair $(a,z)$ for $G$. Since $G_e=\l z\r\cong Z_2$ and $z$ interchanges the vertices $\a$ and $\b$, we see that $G$ acts faithfully on $V\cup E$, and hence $G$ acts faithfully on the underlying graph. Thus our maps satisfy the following hypothesis.

\begin{hypothesis}\label{Hypo}
{\rm
Let $\calM=(V,E,F,\bfI)$ be a locally finite map, with both $|V|\geq3$ and $|F|\geq 3$, so in particular the underlying graph $\Ga=(V,E,\bfI)$ is not $\K_2^{(\lambda)}$ or a simple cycle.
Assume that $\calM$ is $G$-vertex-rotary for some $G\leqslant\Aut\calM$ with a rotary pair $(a,z)$, so that for some arc $(\a,e)$ of $\Gamma$, we have $G=\l a,z\r$,
$G_\a=\l a\r$ with $|a|=k\lambda\geq3$, where $k$ is the valency and $\lambda$ the edge-multiplicity of  $\Ga$. Also $G_e=\l z\r$ with $|z|=2$, $V=[G:\l a\r]$, $E=[G:\l z\r]$ with $e=\l z\r$, $\a=\l a\r$, and $\b=\l a\r z$. Let $f,f'$ be the two faces of  $\calM$ which are incident with the edge $e$. Then  $(\a,e,f)$ and $(\a,e,f')$ are the two flags of $\calM$ incident with the arc $(\a,e)$ and,  replacing $a$ by $a^{-1}$ if necessary, we may assume that $f^a=f'$.
Since $z$ preserves the edge $e$ we also have
\[
z:\ (f,f')\to (f,f')\ \mbox{or}\ (f',f).
\]
}
\end{hypothesis}

Identifying the faces in terms of the group action is the main challenge, and we next  make some comments on our aims for the paper.

\begin{remark}\label{def:rotary-maps}
{\rm
We first give two general constructions of locally finite $G$-vertex-rotary maps for a group $G$ with given rotary pair $(a,z)$: in Construction~$\ref{cons-rotary}$ we construct a $G$-rotary map  $\RM(G,a,z)$ in the case where $|az|$ is finite, and in  Construction~$\ref{cons-bi-rotary}$ we construct a $G$-bi-rotary map $\BRM(G,a,z)$ in the case where $|zz^a|$ is finite.

Locally finite edge-transitive maps are divided into 14 types in \cite{GW,ST}, according to local actions of their automorphism groups; and
vertex-rotary maps form two of these types.
Assume that Hypothesis~\ref{Hypo} holds.
We will show  (in Theorem~\ref{thm}) that the underlying graph $\Gamma$ is isomorphic to $\Cos(G,\l a\r,\l z\r)$, and that $\Gamma$ has two different $G$-vertex-rotary embeddings, depending on the action of $z$ on $(f,f')$.
\begin{itemize}
\item[(i)] If $z$ interchanges $f$ and $f'$, then we show that $|az|$ is finite, the  stabiliser of the face $f$ is $G_f=\l az\r$, and the map is $G$-rotary, sometimes also called an {\it orientably regular} map (see \cite{rotary}). This map is of {\it type $2^P$ex} in the classification scheme of \cite{ST}, and we show that it is isomorphic to $\RM(G,a,z)$.

\item[(ii)]
If  $z$ fixes $f$ and $f'$, then we show that  $\ell:=|zz^a|$ is finite, and the stabilisers of the faces $f$ and $f'$ are $G_f=\l z,z^{a^{-1}}\r$ and $G_{f'}=\l z,z^a\r$. Here $G_f\cong G_{f'}\cong\D_{2\ell}$ and the faces $f, f'$ have equal face-length $2\ell$. This map is {\it $G$-bi-rotary}.  It is of {\it type $2^*$ex} in the scheme \cite{ST}, and we show that it is isomorphic to $\BRM(G,a,z)$.
\end{itemize}
}
\end{remark}

We now give a formal statement of our main result Theorem~\ref{thm} .

\begin{theorem}\label{thm}
Let $\calM=(V,E,F,\bfI)$ be a locally finite $G$-vertex-rotary map such that Hypothesis~$\ref{Hypo}$ holds, and let $f, f'$ be the two faces of $\calM$ incident with the edge $e=\l z\r$,
as in  Notation~$\ref{def:rotary-maps}$. Then either
\begin{enumerate}
\item[(a)] $\calM$ is the $G$-rotary map $\RM(G,a,z)$ of Construction~$\ref{cons-rotary}$ with finite face length $|az|$, and $z$ interchanges $f$ and $f'$; or

\item[(b)] $\calM$ is the $G$-bi-rotary map
$\BRM(G,a,z)$  of Construction~$\ref{cons-bi-rotary}$ with finite face length $2\,|zz^a|$, and $z$ fixes each of $f$ and $f'$.
\end{enumerate}
\end{theorem}

 Theorem~\ref{thm}  follows immediately from Proposition~\ref{p:vrmaps1} proved in Section~\ref{sec:vrmaps}, and we draw attention to the discussion in  Remark~\ref{rem:birot} of bi-rotary maps.  This classification  has some interesting consequences. To start with Propositions~\ref{p:rotary} and~\ref{p:bi-rotary} give the following information about the boundary cycles of vertex-rotary maps.

\begin{corollary}\label{Cor:circular}
Let $\calM=(V,E,F)$ be a locally finite $G$-vertex-rotary map with $|V|\geqslant3$ and
$|F|\geqslant3$. Then the boundary cycles of faces of $\calM$ are all isomorphic regular cycles. In particular, $\calM$ is a circular map if and only if the $G$-action on $V\cup F$ is faithful.
\end{corollary}

The $G$-rotary map  $\RM(G,a,z)$ is defined in Construction~\ref{cons-rotary} as a
coset configuration $\Cos(G,\l a\r,\l z\r,\l az\r)$, and the incidence relation between vertices, edges and faces of the map is given by non-empty intersection of the cosets. In this case the triple of stabilisers $(\l a\r,\l z\r,\l az\r)$, of an incident vertex, edge and face, uniquely determines the family of boundary cycles and hence the rotary map. However this is not the case for the $G$-bi-rotary map $\BRM(G,a,z)$, where the analogous triple of stabilisers $(\l a\r,\l z\r,\l z,z^a\r)$ may correspond to different $G$-bi-rotary maps.

\begin{corollary}\label{cor:iso}
There exist groups $G$ and rotary pairs $(a,z)$ and $(a',z)$ for $G$ such that $(\l a\r,\l z\r,\l z,z^a\r)=(\l a'\r,\l z\r,\l z,z^{a'}\r)$, and $\BRM(G,a,z)\not=\BRM(G,a',z)$.
\end{corollary}

Corollary~\ref{cor:iso} is proved via an explicit example in Proposition~\ref{p:birotary-ex}, see also the discussion in Remark~\ref{rem:bi-rotary}. We mention moreover that the
concepts `rotary' and `bi-rotary' are not intrinsic topological properties of a map. It is possible, for a given map $\calM$, to have different subgroups $X, Y\leq\Aut\calM$ such that $\calM$ is $X$-rotary and $Y$-bi-rotary.
In order to exhibit such maps we use a general group theoretic construction of $G$-flag-regular maps given in Subsection~\ref{s:con-regularmap} in terms of the theory of coset graphs (Section~\ref{def-coset-g}). This allows us to describe explicitly the face boundary cycles for the $G$-flag-regular map, and in particular to identify the map as an $X$-rotary map or $X$-bi-rotary map for a subgroup $X$ of index $2$ in $G$.   We note that  a characterisation of the graphs embeddable as flag-regular maps is available and well known, see \cite[Section 3, especially Theorem 3]{GNSS}. However, for our purposes (for arc-transitive graphs with arbitrary edge-multiplicity), the explicit nature of the group action and the boundary cycles given in Subsection~\ref{s:con-regularmap} is important.

\begin{corollary}\label{rotary=bi-rotary}
There exist a group $X$ with a rotary pair $(a,z)$, and a group $Y$ with a rotary pair $(a',z')$, such that $\RM(X,a,z)=\BRM(Y,a',z')$.
\end{corollary}

We prove Corollary~\ref{rotary=bi-rotary} in Section~\ref{sec:examples} by constructing an infinite family of flag-regular examples with underlying graphs being hypercubes.  See Lemma~\ref{rotary-bi-rotary}.


In Section~\ref{sec:examples}, we construct rotary embeddings and bi-rotary embeddings for hypercubes and complete bipartite graphs providing examples where the face boundary cycles are regular with edge-induced subgraphs $\C_n^{(\lambda)}$ for many positive integers $n$ and $\lambda$. 
Although these concrete examples are all finite maps, we note that there are many infinite
locally finite examples. For example, the universal cover of any finite rotary map of valency $m$ and face-length $n$ with ${1\over 2}+{1\over m}+{1\over n}<1$ is an infinite, locally finite, $G$-rotary map for $G=\l  a, z \mid a^m=z^2=(az)^n =1\r$, see \cite[page\,361]{Conder-survey}.

\section{Coset graphs}\label{def-coset-g}

Let $\Ga=(V,E)$ be a locally finite graph with $|E|>0$.
For a vertex $\alpha\in V$, let $\Gamma(\alpha)$ be the set of vertices adjacent to $\alpha$, sometimes
called the \emph{neighbourhood} of $\a$,  and, as in the introduction, write $E(\alpha) = \{ e\in E\mid \alpha\,\bfI\,e\}$.
If $\Gamma$ has edge-multiplicity $\lambda$, for some positive integer  $\lambda$, and
if the size of $\Gamma(\alpha)$ is independent of  $\alpha$, then we say that $\Gamma$ is \emph{regular},
the constant $|\Gamma(\alpha)|$ is called the \emph{valency} of $\Gamma$, and $|E(\alpha)|=\lambda |\Ga(\a)|$.  For
$e\in E$ incident with vertices $\a$ and $\b$, we sometimes incorporate this information about incidence into our notation and refer to the edge as
$[\alpha,e,\beta]$ or $[\beta,e,\alpha]$. We regard the two arcs $(\alpha,e,\beta)$ and $(\beta,e,\alpha)$ corresponding to $e$ as
`orientations' of $e$. We now define coset graphs with more detail.


\begin{construction}\label{def-m-graph}{\rm
Let $G$ be a group with subgroups  $H$ and $J$ such that $H\ne G$, $|J:H\cap J|=2$, and $|H:H\cap J|$ is finite.
Define an incidence structure $\Cos(G,H,J) = (V, E, \bfI)$ by setting
\[
\begin{array}{lll}
V &=[G:H] =\{Hx\mid x\in G\},\\
E &=[G:J] =\{J x\mid x\in G\},\\
\bfI &= \{(Hx,Jy) \mid yx^{-1}\in JH\} \subseteq V\times E.
\end{array}
\]
We prove that $\Cos(G,H,J)$ is a locally finite graph and call it a {\it coset graph}.
Further, let $g\in J\setminus(H\cap J)$, let $K := H\cap H^g$ so that $H\cap J \leq K \leq H$. Also let  $k := |H:K|$ and
$\lambda := |K:H\cap J|$ (which by assumption are finite), and let $L := \l K,g\r$.
\qed
}
\end{construction}

Our first result Theorem~\ref{thm:cosetgraph} establishes the basic properties of
$\Cos(G,H,J)$. As an easy consequence we prove the first assertion of Theorem~\ref{new-cos-rep} in Corollary~\ref{cor:cosetgraph}. Later in the section we prove Theorem~\ref{thm:cosetgraph2} which allows us to prove Theorem~\ref{cos-rep} and to complete the proof of Theorem~\ref{new-cos-rep}.

\begin{theorem}\label{thm:cosetgraph}
Let $G, H, J, g, K, L, k, \lambda$, be as in Construction~$\ref{def-m-graph}$,
and let $\Gamma$ be the coset graph $\Cos(G,H,J)$.
Then the following statements hold.
\begin{enumerate}
 \item[(a)] $\Gamma$ is a locally finite, regular graph of
valency $k$ and edge-multiplicity $\lambda$; and $\Ga$ is connected if and only if
$G=\l H, J\r$;

\item[(b)] The right multiplication action of $G$ on $V\cup E$ induces an arc-transitive group of automorphisms of $\Gamma$ isomorphic to $G/\Core_G(H\cap J)$; the stabilisers in $G$ of the vertex $H$, edge $J$,
and arc $(H,J,Hg)$, are $H, J$ and $K\cap J=H\cap J$, respectively;

\item[(c)] $L$ satisfies $H\cap L = K$, $|L:K|=2$, and $|L:J|=\lambda$;
Construction~$\ref{def-m-graph}$ applied to $G, H, L$ yields a locally finite simple graph
$\Gamma' := \Cos(G, H, L)$, and moreover, $\Gamma'\cong \SimpCos(G, H, HgH)$;

\item[(d)] $\Gamma'$ is (isomorphic to) the base graph of $\Gamma$;
$\Gamma$ is (isomorphic to) a $(G,\lambda)$-extender of $\Gamma'$; and $\Gamma$ is a simple
graph if and only if $L=J$.

\end{enumerate}
\end{theorem}

Before giving the proof of this result we note the following corollary which proves the first assertion of Theorem~\ref{new-cos-rep}, and we make some remarks about the assumption $\Core_G(H\cap J)$ in that result.

\begin{corollary}\label{cor:cosetgraph}
Using the assumptions and notation of Theorem~$\ref{thm:cosetgraph}$, let $M$ be $\Core_G(H\cap J)$ and write $\overline{x}:=Mx$ for $x\in G$. Then
\begin{enumerate}
\item[(a)] The first assertion of Theorem~$\ref{new-cos-rep}$ is valid.

\item[(b)] The groups $\overline{G}=G/M, \overline{H}=H/M, \overline{J}=J/M$
satisfy the conditions of Construction~$\ref{def-m-graph}$ with $\Core_{\overline{G}}( \overline{H}\cap \overline{J})=1$, and the map $\phi:Hx\to  \overline{H}\overline{x},
\ Jy\to  \overline{J}\overline{y}$, for $x,y\in G$, defines a graph isomorphism from $\Cos(G,H,J)$ to $\Cos(\overline{G}, \overline{H}, \overline{J})$.
\end{enumerate}
\end{corollary}

\proof
(a) If $M=\Core_G(H\cap J)=1$, then the map $\phi$ is just the identity map. It follows from Theorem~\ref{thm:cosetgraph}(b) that $\Cos(G,H,J)$ is $G$-arc-transitive, and from Theorem~\ref{thm:cosetgraph}(a) that $\Cos(G,H,J)$ is non-empty and locally finite.

(b) Let $\Cos(\overline{G}, \overline{H}, \overline{J})= (\overline{V}, \overline{E}, \overline{\bfI})$ be as in Construction~\ref{def-m-graph}. Since $M\subseteq H\cap J$ it is clear that $\phi$ induces bijections $V\to \overline{V}$ and $E\to \overline{E}$. Also since the incidence relation $\overline{\bfI} = \{(\overline{H}\overline{x},\overline{J}\overline{y}) \mid \overline{y}\overline{x}^{-1}\in \overline{J}\overline{H}\}$, it follows that $(Hx, Jy)\in\bfI$ implies that $(\overline{H}\overline{x},\overline{J}\overline{y})\in\overleftarrow{\bfI}$, and hence$\phi$ maps edges of $\cos(G,H,J)$ to edges of   $\Cos(\overline{G}, \overline{H}, \overline{J})$. Thus $\phi$ is a graph isomorphism.
\qed

\begin{remark}\label{rem:m-graph}
{\rm
The reader might wonder why we did not add the requirement `$H\cap J$ core-free' to Construction~\ref{def-m-graph}, since this would remove the possibility that $G$ acts unfaithfully on $\Cos(G,H,J)$. The reason is the importance, in our view, of the observation in Theorem~\ref{thm:cosetgraph}(c). For this, it is essential that we can apply  Construction~\ref{def-m-graph} to the triple $G,H,L$ and construct $\Cos(G,H,L)$.
If we add the requirement `$H\cap J$ core-free' to Construction~\ref{def-m-graph}, then we could only apply this construction to $G, H, L$ if $H\cap L$ is core-free. However it is possible to have $H\cap J$ core-free and $H\cap L$ not core-free.

Here is a simple example. Let $Y=\Sym(\{1,2,3\})$ and $Z\cong \ZZ_\lambda$ for an arbitrary $\lambda>1$, and define $G=Y\times Z, H=Y_1\times Z$, and $J=Y_3\times 1$. Then $G, H, J$ satisfy all the conditions of Construction~\ref{def-m-graph} and in addition $H\cap J=1$ is core-free in $G$. In this case $K=1\times Z$ and $L=J\times Z$ so $\Core_G(L)=Z\ne 1$. We note that $\Cos(G,H,J)\cong \C_3^{(\lambda)}$ and applying Construction~\ref{def-m-graph} to $G,H,L$ (without the core-free requirement) we obtain $\Cos(G,H,L)=\C_3$, the simple base graph.
}
\end{remark}


We prove Theorem~\ref{thm:cosetgraph} in a series of lemmas using the notation introduced in Construction~\ref{def-m-graph} without further reference.
First we collect some technical facts.

\begin{lemma}\label{basic-pty1}
\begin{itemize}
\item[(i)]
 $g^2\in H\cap J\leq K$, and $J=(H\cap J)\l g\r=\l g\r(H\cap J)$with
$H\cap J = K\cap J$ of index $2$ in $J$;
\item[(ii)]  $L = K\l g\r = \l g\r K$ with $K=H\cap L$ of index $2$ in $L$.
Moreover, $J\leq L$ and $\lambda = |L:J|=|K:K\cap J|$;
\item[(iii)] $JH=\l g\r H$, and $(Hx)\,\bfI\,(Jy)$ if and only if $yx^{-1}\in\l g\r H$.
\end{itemize}
\end{lemma}
\proof
Since $|J:H\cap J|=2$ and $g\in J\setminus(H\cap J)$, we have  $g^2\in H\cap J$,  so $J=(H\cap J)\l g\r=\l g\r(H\cap J)$,
with $H\cap J$ normal in $J$ of index $2$. Thus $H^g\cap J = (H\cap J)^g=H\cap J$.
This implies that $K\cap J = H\cap J$, proving part (i). In particular $g\in J\setminus(K\cap J)$ so $g\not\in K$, but also $g^2\in K$
and hence   $L = K\l g\r = \l g\r K$ with $K=H\cap L$ of index $2$ in $L$, proving the first part of (ii).
Now $J=(K\cap J)\l g\r \leq L$, and we have
\[
|L:K\cap J| = |L:K|\cdot |K:K\cap J| = 2 \cdot |K: H\cap J| \ (\mbox{by part (i)})\ = 2\lambda.
\]
Also $|L:K\cap J| = |L:J|\cdot |J:K\cap J| = 2\cdot |L:J|$, so $|L:J|=\lambda$, proving part (ii).
Part (iii) follows from part (i) and the definition of $\bfI$.
\qed

We prove the parts of Theorem~\ref{thm:cosetgraph} in four separate lemmas.

\begin{lemma}\label{coset-graphs}
The claims of Theorem~{\rm\ref{thm:cosetgraph}\,(a)} are valid. Moreover, for  $\Gamma=\Cos(G,H,J)$,
\begin{enumerate}
 \item[(i)]  the edge $Jy$ is incident with the (distinct) vertices $Hy$ and $Hgy$;
\item[(ii)] the vertex $Hx$ is incident with the set $E(Hx)=\{ Jh x\mid h\in H\}$ of $k\lambda$ edges,
and adjacent to the set $\Gamma(Hx) = \{ Hghx \mid h\in H\}$ of $k$ vertices;
\item[(iii)] the set of edges incident with each of two adjacent vertices $Hx, Hghx$
is $\{ Jzhx\mid z\in K\}$, a set of size $\lambda$;
\item[(iv)] the connected component of $\Ga$ containing $Hx$ has vertex set
$\{ Hyx\mid y\in \l H,J\r \}$.
\end{enumerate}
 \end{lemma}

\proof
Since  $yy^{-1}=1\in \l g\r H$ and $y(gy)^{-1}=g^{-1}\in \l g\r H$, and since
$JH=\l g\r H$ (by Lemma~\ref{basic-pty1}), it follows that $Hy\,\bfI\,Jy$ and
$Hgy\,\bfI\,Jy$. Moreover $Hy\ne Hgy$  since $g\not\in H$.

We claim that these are the only vertices incident with $Jy$:  if $Hx\,\bfI\,Jy$,
then $yx^{-1}\in\l g\r H$, and so $x=hg^i y$ for some $h\in H$ and some integer $i$.
Moreover, since $g^2\in H$ (by Lemma~\ref{basic-pty1}), we may take $i\in\{0,1\}$,
and hence  $Hx = Hy$ or $Hgy$ according as $i=0$ or $i=1$ respectively.
This proves the claim, and it follows that  $\Cos(G,H,J)$ is a graph, and part (i) holds.

As we showed above, the edges incident with a given vertex $Hx$ are those cosets $J y$ such that
either $Hx=Hy$ or $Hx = Hgy$, that is to say, either $y=h x$ or $gy=hx$ for some $h\in H$.
Since $g\in J$, we have $Jy=Jgy$, and hence the set
$E(Hx)$ of edges incident with $Hx$ is $\{ Jh x\mid h\in H\}.$
Since for $h, h'\in H$, $Jhx=Jh'x$ if and only if $h'h^{-1}\in H\cap J$, it follows that this set has size
$|H:H\cap J|=k\lambda$.
Further, for $Jhx\in E(Hx)$, the second vertex incident with $Jhx$ is $Hghx$, so
$\Gamma(Hx)=\{ Hghx \mid h\in H\}$. Moreover, $Hghx=Hgh'x$ if and only if $h'h^{-1}\in H\cap H^g = K$, and hence
 $|\Gamma(Hx)|=|H:K|=k$. Since this holds for each $x\in G$, $\Ga$ has valency $k$.
This proves part (ii).

 By part (ii), $Hx\,\bfI\, Jhx$ (for $h\in H$) and by part (i), $Hghx \,\bfI\, Jhx$, so $Hx, Hghx$ are adjacent.
Further, $Hx, Hghx$,
are both incident with an edge $Jh'x$ if and only if $Hghx = Hgh'x$ (by (i)), or
equivalently, $z:=h'h^{-1}\in H\cap H^g = K$. Thus the set of edges
incident with both these vertices is $\{ Jzhx\mid z\in K\}.$   Two edges
$Jzhx, Jz'hx$ in this set are equal if and only if $z'z^{-1}\in K\cap J$,
and so the number of distinct edges incident with both $Hx$ and $Hghx$ is
$|K:K\cap J|$ which is equal to $\lambda$, by Lemma~\ref{basic-pty1}(ii).
This proves part (iii), and also that $\Ga$ has constant edge-multiplicity $\lambda$.
Since by part (ii) $\Ga$ has constant valency $k$, we conclude that $\Ga$ is regular,
Also by parts (ii) and (iii) we have $|E(Hx)|=k\lambda$ for each vertex $Hx$, so $\Ga$ is locally finite.

Noting that $JH=\l g\r H$, by Lemma~\ref{basic-pty1}(iii), it follows from part (ii) that $Hz$
lies in the connected component of $\Ga$ containing $Hx$ if and only if $Hz=Hyx$ for some
$y\in\l H,J\r$. Thus part (iv) is proved, and $\Ga$ is connected if and only if  $\l H,J\r=G$.
Therefore all the claims of Theorem~\ref{thm:cosetgraph}\,(a) are proved.
\qed

\begin{lemma}\label{coset-graphs-2}
The claims of Theorem~{\rm \ref{thm:cosetgraph}\,(b)} are valid. Moreover
the kernel of the $G$-action on the edges of $\Cos(G,H,J)$ is $\Core_G(J)$, and either
$\Core_G(J)$ also acts trivially on vertices, or there exists $\lambda$ such that
each connected component of $\Ga$ is the $\lambda$-extender $\K_2^{(\lambda)}$
of a complete graph on two vertices.
\end{lemma}

\proof
We claim that $G$ acting by right multiplication induces a subgroup of automorphisms of $\Gamma =\Cos(G,H,J)$.
Let $z\in G$. Then, for $x, y\in G$,
\[
\begin{array}{llll}
Hx\ \bfI\ Jy & \Longleftrightarrow & yx^{-1}\in \l g\r H &\mbox{(by Lemma~\ref{basic-pty1}(iii))}\\
 &\Longleftrightarrow & (yz)(xz)^{-1}\in \l g\r H \\
 &\Longleftrightarrow & (Hx)z\ \bfI \ (J y)z.
\end{array}
\]
Thus right multiplication by $z$ preserves $\bfI$, proving the claim. The stabiliser of the vertex $H$, edge $J$, and arc $(H,J,Hg)$ under this action
are $H$, $J$, and $H\cap J\cap H^g = K\cap J$, respectively. Now $K\cap J\leq H\cap J$ since $K\leq H$. Conversely $H\cap J$ fixes both of the vertices $H$ and $H^g$ incident with the edge $J$ and the reverse inclusion holds, so $K\cap J=H\cap J$. The kernels of the $G$-actions
by right multiplication on $V$ and on $E$ are $\Core_G(H)$ and $\Core_G(J)$, respectively. Hence
the kernel of the $G$-action on $V\cup E$ is $\Core_G(H)\cap \Core_G(J)$, and it is straightforward to check that this subgroup is equal to $\Core_G(H\cap J)$.

To prove the claims of Theorem~{\rm \ref{thm:cosetgraph}\,(b)} it remains to show
that $G$ acts arc-transitively. By Lemma~\ref{coset-graphs}~(i),
each arc of $\Gamma$ is of the form $(Hy, Jy, Hgy)$ (for, since $g\in J$ and $g^2\in H\cap J$, the arc
$(Hgy, Jy, Hy)$ is equal to $(Hy', Jy', Hgy')$ with $y'=gy$). Right
multiplication of the arc $(Hy, Jy, Hgy)$ by $y^{-1}$ maps it to $(H,J,Hg)$.  Thus $G$ has just one orbit on arcs.

As noted above, the kernel of the $G$-action on edges is
$\Core_G(J)$. Suppose that $\Core_G(J)$ acts
nontrivially on the vertices of $\Gamma$. Then since $G$ is vertex-transitive, all
$\Core_G(J)$-orbits on vertices have the same length, say $a>1$.
By Construction~\ref{def-m-graph}, $H$ is a proper subgroup of $G$ so
$\Ga$ has at least two vertices and, by Lemma~\ref{coset-graphs}\,(i),
for each $y\in G$, $e_y:= [Hy, Jy, Hgy]$ is an edge.  However $e_y$ is fixed by
$\Core_G(J)$, and so $\Core_G(J)$ either fixes on interchanges the two vertices
$Hy, Hgy$ incident with $e_y$. Since each $\Core_G(J)$-orbit on vertices has size $a>1$, it follows that $a=2$ and $\Core_G(J)$ interchanges $Hy, Hgy$ for each $y\in G$.
By Theorem~\ref{thm:cosetgraph}(a) (which holds by Lemma~\ref{coset-graphs}),
$\Gamma$ is regular of valency $k$. If $k>1$ then there are edges $e_1 = [H,J,Hg],
e_h = [H,Jh, Hgh]$ incident with $H$ and with $Hg\ne Hgh$, for some $h\in H$, by Lemma~\ref{coset-graphs}.
We have shown that $\Core_G(J)$ fixes setwise $\{H, Hg\}$ and also $\{H, Hgh\}$, and interchanges the two vertices  in each of these sets, which is impossible. Hence
$\Ga$ has valency $k=1$. Thus, since $G$ is arc-transitive on $\Gamma$, there is a constant $\lambda$
such that each connected component of $\Ga$ is isomorphic to $\K_2^{(\lambda)}$.
\qed

\begin{remark}\label{rem:core}
{\rm
By Lemma~\ref{coset-graphs-2}, if $\Ga=\Cos(G,H,J)$ is connected and has at least three vertices,
then the group of automorphisms of $\Ga$ induced by $G$ is isomorphic to the group
that $G$ induces on the edge-set of $\Ga$.
In contrast to this,  there are connected coset graphs $\Cos(G,H,J)$ with arbitrarily large vertex sets, and
arbitrarily large kernels of the $G$-action on vertices. For example, for any
$n$ and $\lambda$, the $\lambda$-extender $\Ga=\C_n^{(\lambda)}$ of a simple cycle $\C_n$ of length $n$
is a coset graph $\Cos(G,H,J)$, where $G =
\D_{2n}\times\Sym(\lambda)$, $H=\ZZ_2\times\Sym(\lambda)$, $J=\l g\r\times \Sym(\lambda-1)$
for some $g\in \D_{2n}\setminus H$, and the kernel of the vertex action is $\Core_G(H)= \Sym(\lambda)$.
}
\end{remark}

\begin{lemma}\label{coset-graphs-3}
The claims of Theorem~{\rm \ref{thm:cosetgraph}\,(c)} are valid.
\end{lemma}

\proof
We have $\Ga = \Cos(G, H, J) = (V, E, \bfI)$.
The properties $H\cap L=K$ and $|L:K|=2$ follow from Lemma~\ref{basic-pty1}~(ii).
Moreover, this result also gives $L=K\l g\r$, so $g\in L\setminus (H\cap L)$.
Thus we may apply Construction~\ref{def-m-graph}, with $J$ replaced by $L$,
to obtain $\Gamma' =\Cos(G, H, L) = (V', E', \bfI')$. Note that in this construction $V'=V$ and we may choose the same
element $g$ as in the construction of $\Ga$, and hence we have the same subgroup
$K = H\cap H^g$. Also since $L=\l K, g\r$, the subgroup $L$ has the same properties for $\Ga'$
as it has for $\Ga$. In particular it follows from  Lemma~\ref{coset-graphs}, applied to $\Ga'$,
that $\Ga'$ is locally finite with constant valency $k$ and constant edge-multiplicity $|L:L|=1$, so $\Gamma'$
is a regular simple graph.

By Lemma~\ref{coset-graphs-2}, applied to $\Ga'$, $\Gamma'$ is $G$-arc-transitive and the element
$g\in G$ maps the arc $(H,L,Hg)$ of $\Ga'$ to its reverse arc $(Hg,L,H)$ (since $g\in L$ and $g^2\in H$).
Thus all the conditions hold for constructing  the simple coset graph $\SimpCos(G,H,HgH) = (V'', E'', \bfI'')$.
Recall that $V''=V=V'$, and two vertices $Hx, Hy$ are adjacent in $\SimpCos(G,H,HgH)$ if and only if $yx^{-1}\in HgH$,
that is, $y=h'ghx$ for some elements $h',h\in H$.
Thus the edge set $E''$ consists of all vertex pairs of the form
\[
\{Hx,Hghx\},\ \mbox{for some $x\in G$ and $h\in H$,}
\]
and incidence $\bfI''$ between vertices and edges is given by inclusion.

By Lemma~\ref{coset-graphs}, applied to $\Ga'$ and recalling that $g\in L$, both of the vertices $Hx$ and $Hghx$ are
incident in $\Gamma'$ with the edge $Lghx$, and each edge of
$\Gamma'$ is of the form $[Hx, Lghx, Hghx]$ for some $x\in G$, $h\in H$.
Define $f: V\cup E'' \rightarrow V\cup E'$ such that
$f$ restricts to the identity map on $V$, and $f$ sends the edge $\{Hx,Hghx\}$ of $E''$
to the edge $Lghx$ of $\Gamma'$. It follows from our comments that
$f$ is a bijection and induces a bijection from  $\bfI''$ to $\bfI'$. Hence $f$ is a graph isomorphism from  $\SimpCos(G,H,HgH)$ to $\Gamma'$.
\qed

\begin{lemma}\label{coset-graphs-4}
The claims of Theorem~{\rm\ref{thm:cosetgraph}\,(d)} are valid.
\end{lemma}

\proof
It follows from the definition of a base graph,
and the proof of Lemma~\ref{coset-graphs-3}, that the base graph of $\Ga$ is $\Ga'=\Cos(G,H,L)
\cong \SimpCos(G,H,HgH)$. Also, by Lemma~\ref{coset-graphs}\,(iii), $\Gamma$ has constant
edge-multiplicity $\lambda$ and hence $\Gamma$ is isomorphic to a $(G,\lambda)$-extender of $\Gamma'$.
Finally, $\Gamma$ is a simple graph if and only if $\lambda=1$ and,
since $\lambda = |L:J|$  by Lemma~\ref{coset-graphs},
this holds if and only if $L=J$.
\qed

\noindent
\emph{Proof of Theorem~{\rm \ref{thm:cosetgraph}}.}\quad This follows from Lemmas~\ref{coset-graphs},
\ref{coset-graphs-2}, \ref{coset-graphs-3}, and \ref{coset-graphs-4}.
\qed


Next we demonstrate the universality of Construction~\ref{def-m-graph}
by showing that every arc-transitive graph (with arbitrary finite edge-multiplicity)
can be represented as a coset graph arising from Construction~\ref{def-m-graph}.

\begin{theorem}\label{thm:cosetgraph2}
Let $\Ga=(V,E,\bfI)$ be a locally finite graph with $|E|>0$, and suppose
that $G\leqslant\Aut\Ga$ is such that $\Ga$ is $G$-arc-transitive with valency $k$
and edge-multiplicity $\lambda$. Let $\a, \b$ be a pair
of adjacent vertices, let $e$ be an edge incident with $\a$ and $\b$,
and let $H=G_\a$ and $J=G_e$.
\begin{enumerate}
\item[(a)] Then $|J:H\cap J|=2$, $|H:H\cap J|$ is finite,
and  $\Ga\cong\Cos(G,H,J)$;
\item[(b)] and if  $g\in G$ reverses the arc $(\a,e,\b)$ then,
in the notation of Construction~$\ref{def-m-graph}$,
\[
G_{\a\b} = K := H\cap H^g,\quad G_{\{\a,\b\}} = L := K\l g\r,\quad G_e = J = (K\cap J)\l g\r,
\]
$\lambda = |L:J| = |K:(K\cap J)|$, and $k=|H:K|$.
\end{enumerate}
\end{theorem}

\proof
First we establish the claims about subgroup equalities and subgroup indices.
Since the element $g$ reverses the arc $(\a,e,\b)$, we have $\a^g=\b$ and hence
$K = H\cap H^g = G_\a\cap G_\b = G_{\a\b}$, so $k=|H:K|$, and $L = K\l g\r = G_{\{\a,\b\}}$.
Moreover, since $G$ is arc-transitive, it follows that $K = G_{\a\b}$ is transitive
on the $\lambda$ edges incident with $\a$ and $\b$, and its index $\lambda$
subgroup stabilising $e$ is $G_{\a\b e}=K\cap J$. In addition, $g$ fixes the edge $e$,
so $g\in J\setminus (K\cap J)$, $G_e = J = (K\cap J)\l g\r$, and we have $|L:K|=|J:K\cap J|=2$.
Since $\Ga$ is locally finite, the subgroup $G_{\a e}=H\cap J$ has finite index in $H=G_\a$, and
since  $e$ is incident with precisely the two vertices $\a, \b$, the group
$G_{\a e}$ also fixes $\b$ and it follows that $K\cap J=H\cap J$.
Thus $|J:H\cap J|=2$. Finally $|L:J| = |L:K|\cdot |K:K\cap J|/|J:K\cap J| = |K:K\cap J|$, which
equals $\lambda$.

Thus $G, H, J$ satisfy the
conditions of Construction~\ref{def-m-graph}, and we obtain a coset graph $\Ga' := \Cos(G, H, J) = (V', E', \bfI')$,
where $V' = [G:H], E' = [G:J],$ and $(Hx)\,\bfI'\,(Jy)$ if and only if $yx^{-1}\in JH$.
By Lemma~\ref{basic-pty1}, $JH=\l g\r H = H\cup gH = H \cup g^{-1}H$ (since $g^2\in H\cap J$).
We will prove that the following map defines a graph isomorphism from $\Ga$ to $\Ga'$:
\[
f:\ \ \begin{array}{l}
\a^x\mapsto H x,\\
e^x\mapsto Jx,
\end{array}\ \ \mbox{for $x\in G$,}
\]
that is to say, we will prove that $f$
is a well defined bijection $f: V\cup E \rightarrow V'\cup E'$ such that $(V)f=V'$, $(E)f=E'$,
and $\a^x \ \bfI\ e^y$ if and only if $(\a^x)f\ \bfI'\ (e^y)f$.
First we note that
\[\begin{array}{lll}
\a^x= \a^y  	& \Longleftrightarrow & \a^{xy^{-1}}=\a \\
		& \Longleftrightarrow & xy^{-1}\in G_\a=H\\
		& \Longleftrightarrow & Hxy^{-1} =H\\
		& \Longleftrightarrow & Hx = Hy.
\end{array}\]
From this, and the fact that each coset $Hx$ occurs in the image of $f$,
it follows that the restriction $f|_V$ is a well defined bijection $V\rightarrow V'$.
The same argument, interchanging $V, V'$ with $E, E'$, shows that
$f|_E$ is a well defined bijection $E\rightarrow E'$. Thus $f$ is a bijection and
$(V)f=V'$, $(E)f=E'$. Now we check the assertion about incidence. Let $\a^x\in V$ and $e^y\in E$.
Then
\[\begin{array}{lll}
\a^x\ \bfI\ e^y   & \Longleftrightarrow & \a^{xy^{-1}}\ \bfI\ e \\
		  & \Longleftrightarrow & \a^{xy^{-1}}=\a\ \mbox{or}\ \a^{xy^{-1}}=\b = \a^g\\
		  & \Longleftrightarrow & xy^{-1} \in H\  \mbox{or}\ xy^{-1}\in Hg\\
		  & \Longleftrightarrow & xy^{-1} \in H\cup Hg \\
		  & \Longleftrightarrow & yx^{-1} \in H\cup g^{-1}H = JH\\
		  & \Longleftrightarrow & Hx\ \bfI'\ Jy.\\
\end{array}\]
Thus $\Ga\cong\Ga'$, completing the proof.
\qed

We now complete the proof of Theorem~\ref{new-cos-rep}, and prove Theorem~\ref{cos-rep}.
\vskip0.1in

\noindent
\emph{Proof of Theorem~\ref{new-cos-rep}.}\quad The first assertion of Theorem~\ref{new-cos-rep} follows from Corollary~\ref{cor:cosetgraph}(a), while the converse assertion is proved in Theorem~\ref{thm:cosetgraph2}(a).
\qed

\vskip0.1in
\noindent
\emph{Proof of Theorem~\ref{cos-rep}.}\quad Let $\Gamma=\Cos(G,H,J)$ for some $G, H, J$ as in Construction~\ref{def-m-graph}, with $H\cap J$ core-free. Then by Theorem~\ref{thm:cosetgraph}(a) and (b), $\Gamma$ is non-empty, locally finite, and $G$-arc-transitive. The fact that $J\leq L$ follows from Lemma~\ref{basic-pty1}(ii), and the other assertions follow from Theorem~\ref{thm:cosetgraph}(c) and (d).
\qed

In Theorem~\ref{thm:cosetgraph}(c), we saw that, for $\Ga = \Cos(G,H,J)$ as in
Construction~\ref{def-m-graph}, we could also apply the construction with $J$ replaced by $L$, obtaining the simple base graph of $\Ga$.
Finally in this section we explore graphs obtained from $\Ga$ by replacing $J$ with various other subgroups. Recall the definition of a $(G,\mu)$-extender of a graph.

\begin{lemma}\label{G-extender}
Let  $G, H, J, g, K, L, \lambda$ be as in Construction~\ref{def-m-graph},
let  $\Ga=\Cos(G,H,J)$, and let $\Sigma=\Cos(G,H,L)$ denote the simple base graph of $\Ga$ (see Theorem~\ref{cos-rep}(a)).
\begin{itemize}
\item[(i)] If, for some odd integer $b$, $\l g^b\r\leqslant J'\leqslant J$ and $\mu=|J:J'|$ is finite, then
$|J': H\cap J'|=2, |H:H\cap J'|$ is finite,   $\Ga':= \Cos(G,H,J')$
is $G$-arc-transitive,  $\Ga'$ is a $(G,\mu)$-extender of $\Ga$, and $\Ga'$ is a $(G,\lambda\mu)$-extender of
$\Sigma$.

\item[(ii)] If $R$ is a proper subgroup of $K$ with $\mu = |K:R|$ finite, and such that, for some odd integer $b$, $g^{2b}\in R$ and $R^{g^b}=R$, then  $\Cos(G, H, R\l g^b\r)$ is a $(G,\lambda\mu)$-extender of the simple graph $\Sigma$. Moreover, if $K$ is finite, then $\Sigma$ has a non-trivial $G$-extender unless $L=\l g\r\cong Z_{2^a}$ for some $a\geq1$.

\end{itemize}
\end{lemma}
\proof
The graph $\Ga$ is independent of the choice of the element $g\in J\setminus(H\cap J)$
(by Lemma~\ref{basic-pty1}).

(i) Since $g^b\in J'$ for some odd integer $b$, we may replace $g$ by $g^b$
in Construction~\ref{def-m-graph} so that, by Lemma~\ref{basic-pty1}(i), the product $J'(H\cap J)=J$, and hence $|J':(H\cap J')|=2$ and $|H\cap J:H\cap J'|=|J:J'|=\mu$. Thus $|H:H\cap J'|=\mu |H:H\cap J|$
is finite, so  $G, H, J'$ satisfy the conditions of Construction~\ref{def-m-graph}. (Note also that, since $g^2\in H$, we have $K=H\cap H^g = H\cap H^{g^b}$.)  Therefore,
by Theorem~\ref{thm:cosetgraph}(a), $\Ga' := \Cos(G, H, J')$ is $G$-arc transitive
with edge-multiplicity $|L:J'|=\lambda\mu$. It follows from Theorem~\ref{thm:cosetgraph}(d)
that $\Ga'$ is a $(G,\mu)$-extender of $\Ga$, and a  $(G,\lambda\mu)$-extender of $\Sigma$.

(ii)
Suppose that $R < K$, $\mu:=|K:R|$
is finite, $g^{2b}\in R$ and $R$ is normalised by $g^b$, for some odd integer $b$.
Since $b$ is odd, $g^b\not\in K$ by Lemma~\ref{basic-pty1}, and hence $g^b\not\in R$. Then since $g^{2b}\in R$ and $R^{g^b}=R$, the subgroup $J':= R\l g^b\r$ satisfies $|J':R|=2$. Further, since $b$ is odd, it follows from Lemma~\ref{basic-pty1}(ii) that $J=K\l g^b\r$.  Thus
$|J:J'|=|K\l g^b\r : R\l g^b\r| = |K:R|=\mu >1$, and so $J'$ satisfies the conditions of part (i).
Hence by part (i), $\Cos(G,H,J')$ is a $(G,\lambda\mu)$-extender of $\Sigma$.

Finally suppose that $K$ is finite. By Theorem~\ref{thm:cosetgraph}, $\Ga$ is a $(G,\lambda)$-extender of $\Sigma$, and this is a non-trivial $G$-extender of $\Sigma$ if $\lambda>1$. So suppose that $\lambda=1$. Then $L=J$ by Lemma~\ref{basic-pty1}. If $K$ contains $\l g^2\r$ as a proper subgroup then $\l g^2\r$ satisfies all the conditions for $R$ in the previous paragraph, taking $b=1$, and this yields a non-trivial $G$-extender of $\Sigma$. Thus we may assume further that $K=\l g^2\r$ and hence that $L=J=\l g \r$.  Since $K$ is finite, also $|g|$ is finite. If  $|g|$ is divisible by an odd integer $b>1$ then the subgroup $R:=\l g^{2b}\r$
satisfies all the condition of the previous paragraph for the integer $b$, again giving a non-trivial $G$-extender of $\Sigma$. This leaves the exceptional case where $L=\l g\r$ is cyclic of $2$-power order.
\qed

As an example we apply this theory to the Petersen graph, and demonstrate that the number of non-trivial $G$-extenders of a $G$-arc transitive graph may depend on the choice of the arc-transitive subgroup $G$.

\begin{example}\label{Petersen-1}
{\rm
Let $\Ome=\{1,2,3,4,5\}$, and let $\Sig$ be the Petersen graph.
Then the ten vertices of $\Sig$ can be identified with the ten 2-subsets of $\Ome$ such that two vertices are adjacent
if and only if the corresponding 2-subsets are disjoint.
 Consider the adjacent vertices $\a=\{2,3\}$ and $\b=\{4,5\}$.

(1). Let $G=\Alt(\Ome)=\A_5$.
Then the element $g=(24)(35)$ interchanges $\a$ and $\b$, and the stabilisers $H=G_\a=\l (145),(23)(45)\r=\D_6$,
 $K=G_{\a\b}=\l (23)(45)\r=\ZZ_2$, and $L=G_{\{\a,\b\}}=\l b,g\r=\ZZ_2^2$, where $b=(23)(45)$.
Hence, by Lemma~\ref{coset-graphs-2}, $\Sig\cong \Cos(G,H,L)=\Cos(\A_5,\D_6,\D_4)$.
There are exactly two elements which interchange $\a,\b$, namely $g=(24)(35)$ and $bg=(25)(34)$.
Thus  $\Sig$ has 2 different $\A_5$-arc-transitive 2-extenders: namely
$\Cos(G,H,\l g\r)$ and $\Cos(G,H,\l bg\r)$.

(2). Alternatively, let $G=\Sym(\Ome)=\S_5$.
Then $G=\Aut\Sig$, and the stabilisers are $H=G_\a=\l (145),(23), (45)\r=\D_{12}$, and $L=G_{\{\a,\b\}}=\l (23),(45),g\r=\D_8$.
Hence also $\Sig\cong \Cos(G,H,L)=\Cos(\S_5,\D_{12},\D_8)$.
Here there are four elements in $G$ that interchange $\a$ and $\b$, namely $(24)(35)$, $(25)(34)$, $(2435)$ and $(2534)$.
Applying Lemma~\ref{G-extender}(i), we see that $\Sig$ has two different $\S_5$-arc-transitive 4-extenders:
$\Cos(G,H,\l (24)(35)\r)$ and $\Cos(G,H,\l (25)(34)\r)$.
\qed
}
\end{example}

\section{Vertex-rotary graphs}\label{s:vertrot}

In this section we study vertex-rotary graphs. Our first result Lemma~\ref{lem:rotpair} establishes in particular that each such graph corresponds to a rotary pair as in Definition~\ref{def-rotarypair}.

\begin{lemma}\label{lem:rotpair}
Let $\Ga=(V,E,\bfI)$ be a connected locally finite graph with $E\ne\emptyset$, and suppose
that $G\leqslant\Aut\Ga$ is such that $\Ga$ is $G$-vertex-rotary with valency $k$
and edge-multiplicity $\lambda$. Let $\a, \b$ be a pair
of adjacent vertices, let $e$ be an edge incident with $\a$ and $\b$,
and let $H=G_\a$ and $J=G_e$. Then
\begin{enumerate}
\item[(a)] $\Ga\cong\Cos(G,H,J)$ and $G=\l a,z\r$ where $H=\l a\r\cong Z_{k\lambda}$, $H\cap H^z=\l a^k\r\cong Z_\lambda$, $J=\l z\r\cong Z_2$, and $H\cap J=1$;
\item[(b)] $z, z^a\not\in\l a\r$, $(a,z)$ is a rotary pair for $G$, and $G$ is regular on the arc set of $\Ga$.
\end{enumerate}
\end{lemma}

\proof
By definition a $G$-rotary graph is $G$-arc-transitive, and hence by
Theorem~\ref{thm:cosetgraph2}, we may assume that $\Ga=\Cos(G,H,J)$, and we have $|J:H\cap J|=2$. Also
by Theorem~\ref{thm:cosetgraph}(a), $G=\l H,J\r$ since $\Ga$ is connected. It follows from the definition of a $G$-vertex-rotary graph that $H$ induces a regular cyclic group on $E(\a)$ so $H\cap J$ is normal in $H$, $H\cap J$ fixes $E(\a)$ pointwise, and $H/(H\cap J)\cong Z_{k\lambda}$. In particular there exists $a\in H$ such that $H=\l a,H\cap J\r$.

We claim that $H\cap J=1$. Now for each $\b$ adjacent to $\a$, and each edge $e\in E(\a)$ incident with $\b$, $H\cap J$ fixes $e$ and hence fixes $\b$. Then as $G_\b$ acts regularly on $E(\b)$ and $H\cap J\leq G_\b$ fixes $e\in E(\b)$, it follows that $H\cap J$ fixes $E(\b)$ pointwise also. Since $\Ga$ is collected we conclude that $H\cap J$ fixes each vertex and each edge of $\Ga$ and hence $H\cap J=1$, proving the claim.  Thus
$J=\l z\r\cong Z_2$, $H=\l a,H\cap J\r=\l a\r\cong Z_{k\lambda}$, and $G=\l H,J\r=\l a,z\r$. As $\a^z=\b$ it follows that $G_{\a\b}=H\cap H^z$ has index $k=|\Gamma(\a)|$ in $H$ and hence $H\cap H^z=\l a^k\r\cong Z_\lambda$, and part (a) is proved.

By definition $G$ is arc-transitive on $\Ga$ and as $H\cap J$ is the stabiliser of the arc $(\a,e)$ and $H\cap J=1$, it follows that $G$ is arc-regular. Now $z\not\in\l a\r$ since $H\cap J=1$, and hence $(a,z)$ is a rotary pair. Also, if $z^a\in\l a\r$ then $\l a\r$ would also contain $az^a a^{-1}=z$, which is a contradiction. This  proves part (b).
\qed

The following converse to Lemma~\ref{lem:rotpair} establishes the equivalence between vertex-rotary graphs and rotary pairs.

\begin{lemma}\label{lem:rotpair2}
Let $G=\l a,z\r$ be a group with rotary pair $(a,z)$ as in Definition~$\ref{def-rotarypair}$, and let $H=\l a\r$ and $J=\l z\r\cong Z_2$. Then  $\Ga:=\Cos(G,H,J)$ is a $G$-vertex-rotary graph of valency $k:=|H:H\cap H^z|$ and edge-multiplicity $\lambda:=|H\cap H^z|$, and $|a|=k\lambda$. Further, if $|a|\geq3$ and $H\ne H^z$, then $\Gamma$ is not $\K_2^{(\lambda)}$ and $\Gamma$ is not a simple cycle.
\end{lemma}

\proof
By Definition~$\ref{def-rotarypair}$, $|a|$ is finite, and $z\not\in H$ so $H\cap J=1$. Thus $G, H, J$ satisfy the conditions of Construction~\ref{def-m-graph} with $H\cap J$ core-free in $G$. Hence by Theorem~\ref{thm:cosetgraph} (b), $\Ga =\Cos(G,H,J)$ is $G$-arc-transitive. Also, by Lemma~\ref{coset-graphs}(ii), for $\alpha=\l a\r$, $E(\alpha)=\{Ja^i\mid 0\leq i<|a|\}$. Moreover the cosets $Ja^i$ in $E(\alpha)$ are pairwise distinct since $z\not\in\l a\r$, and hence $|E(\alpha)|=|a|$. Thus $G_\alpha=H$ acts transitively and cyclically on $E(\alpha)$, and hence $\Gamma$ is $G$-vertex-rotary. By Lemma~\ref{lem:rotpair}(a), $|a|=k\lambda$,  $k=|H:H\cap H^z|$ and $\lambda:=|H\cap H^z|$.
Finally, if $\Gamma=\K_2^{(\lambda)}$ then $H$ fixes both vertices so $H=H^z$, while if $\Gamma$ is a simple cycle then $k=2$ and $\lambda=1$ so $|a|=2$.
\qed

If the group $G$ in Lemma~\ref{lem:rotpair} is abelian, or if the parameters $k, \lambda$ are very small, then the graph can be easily identified.

\begin{corollary}\label{cor:rotpair}
Let $G=\l a,z\r$ and $\Gamma$ be as in Lemma~$\ref{lem:rotpair}$.
\begin{enumerate}
\item[(a)] If $G$ is abelian then $\Ga\cong \K_2^{(\lambda)}$ and $G=\l a\r\times \l z\r\cong Z_\lambda\times Z_2$ with $k=1$;
\item[(b)] and conversely if  $\Ga\cong \K_2^{(\lambda)}$ (so $k=1$), then   $G=\l a\r:\l z\r$, and $a^z=a^j$ for some integer $j$ satisfying $1\leq j\leq \lambda$, $\gcd(j,\lambda)=1$, and $j^2\equiv 1\pmod{\lambda}$;
\item[(c)] if $k\lambda=2$ (so $a$ and $z$ are both involutions), then either
	\begin{enumerate}
	\item[(i)] $k=1, \lambda=2$, $\Gamma=\K_2^{(2)}$, and $G=\l a\r\times \l z\r\cong Z_2\times Z_2$; or
	\item[(ii)] $k=2, \lambda=1$, $\Gamma=\C_r$, for some $r\geq3$, and $G=\l a, z\r\cong\D_{2r}$ with $|az|=r$.
\end{enumerate}	
\end{enumerate}
\end{corollary}

\proof
Suppose that $G$ is abelian, that is, $az=za$, and $G=\l a\r\times \l z\r\cong Z_{k\lambda}\times Z_2$ (since $z\not\in\l a\r$ by Lemma~\ref{lem:rotpair}(b)).  This implies that $|V|=|G:H|=2$ and it follows that $k=1$ and $\Ga=\K_2^{(\lambda)}$, proving part (a). Conversely suppose that $k=1$ and $\Ga\cong \K_2^{(\lambda)}$. Then $V=\{\a,\b\}$, $E=E(\a)$ and $G\leq \Aut(\Ga)=\Sym(E)\times\Sym(V)=\Sym(\lambda)\times\Sym(2)$.   Now $|G:H|=|V|=2$ so $H$ is normal in $G$ and hence $a^z\in H=\l a\r\cong Z_{\lambda}$,  so $a^z=a^j$ for some integer $j$ satisfying $1\leq j\leq \lambda$, and $\gcd(j,\lambda)=1$. Since $z^2=1$, also $j^2\equiv 1\pmod{\lambda}$, and part (b) is proved.

Finally suppose that $k\lambda=2$, so $|a|=|z|=2$ and hence $G\cong \D_{2r}$ where $r=|az|$, noting that $r\geq2$ by Lemma~\ref{lem:rotpair}(b). Thus $|V|= |G:H|=|G:Z|=|E|=r$. If $k=1$ then $\lambda = k\lambda=2$ and as $\Gamma$ is connected,  $|V|=r=2$. Then, as in the previous paragraph we conclude that $\Ga=\K_2^{(\lambda)}$, and  part (c)(i) holds. So we may assume that $k=2$ and $\lambda=1$. This implies that $|V|=r\geq k+1=3$, and that $\Gamma$ is a simple cycle $\C_r$, and part (c)(ii) holds.
\qed

\begin{remark}\label{rem:k2lambda}
{\rm
We make some remarks about the situation in case (b) of Corollary~\ref{cor:rotpair}, where it is possible for $G$ to be nonabelian. Let $V=\{\a,\b\}$ and $\tau=(\a,\b)$,
the generator of $\Sym(V)$. Since $a$ fixes both vertices, it has the form $a=(a_0,1)\in \Sym(E)\times\Sym(V)$ with $a_0$ a $\lambda$-cycle in $\Sym(E)$, and since $z$ interchanges $\a$ and $\b$, $z=(z_0,\tau)$ for some $z_0\in\Sym(E)$ and  $z_0^2=1$ since $z^2=1$. If $j=1$ then $z_0$ centralises the $\lambda$-cycle $a_0$ and fixes $e\in E$, and this implies that $z_0=1$ so $G$ is abelian, and in fact $G$ is cyclic if $\lambda$ is odd. On the other hand if $j=\lambda-1$ then $z_0$ inverts $H$ and so $G$ is a dihedral group of order $2\lambda$. There are other possibilities: for example if $\lambda=8$ and $a_0=(12345678)$, then we could have $j=3$ with $z_0=(24)(37)(68)$, or
$j=5$ with $z_0=(28)(37)(46)$. These extra cases only arise for composite values of $\lambda$.
}
\end{remark}

By Lemma~\ref{lem:rotpair}, we may identify a $G$-vertex-rotary graph $\Ga=(V,E,\bfI)$
with the connected, $G$-arc-regular coset graph $\Ga=\Cos(G,\l a\r, \l z\r)$, where $(a,z)$ is a rotary pair for $G$. Here  $V=[G:\l a\r]$ is the vertex set and $E=[G:\l z\r]$ is the edge set of $\Ga$, $\a$ is the vertex $\l a\r$, $\beta$ is $\l a\r z$, and the edge $[\a,e,\b]$ is $\l z\r$. Moreover we assume that $\Ga$ has valency $k$ and edge-multiplicity $\lambda$, so $|a|=k\lambda$. We now derive more information about the subgroup structure of $G$. For a group $G$ acting on a set $\Delta$, the \emph{kernel} $G_{(\Delta)}$ is the subgroup of $G$ fixing each point of $\Delta$.

\begin{lemma}\label{kernel-V-F}
With the notation in the previous paragraph, the following hold.
\begin{enumerate}
\item[(a)] The kernel of the $G$-action on vertices is $G_{(V)}=\l a\r\cap\l a^z\r$.

\item[(b)] The intersections $\l a\r\cap\l az\r$ and $\l a\r\cap\l zz^a\r$ are normal subgroups of $G$.

\item[(c)] If either $\l a\r\cap\l az\r$ has index at most $2$ in $\l az\r$, or $\l a\r\cap\l z,z^a\r$ has index at most $2$ in $\l z,z^a\r$,
then $G=\l a\r{:}\l z\r$ and $\Ga=\K_2^{(\lambda)}$ with $\lambda=|a|$.
\end{enumerate}

\end{lemma}

\proof
(a) Now $\b=\l a\r z=\a^z$, a vertex adjacent to $\a$, so $G_\b=G_\a^z=\l a^z\r$.
The kernel $G_{(V)}$ fixes both $\a$ and $\b$, and so $G_{(V)}\unlhd G_\a\cap G_\b=\l a\r\cap\l a^z\r$.
As $\l a\r$ is abelian, its subgroup $\l a\r\cap\l a^z\r$ is normalised by $a$, and also,  clearly, $\l a\r\cap\l a^z\r$ is normalised by $z$. Thus $\l a\r\cap\l a^z\r$ is a normal subgroup of $\l a,z\r=G$. Since $\l a\r\cap\l a^z\r$ fixes $\a$, it follows that
$\l a\r\cap\l a^z\r$ fixes all vertices in $V$.
Therefore, $\l a\r\cap\l a^z\r=G_{(V)}$.

(b)
The intersection $\l a\r\cap\l az\r$ is contained in the cyclic subgroups $\l a\r$ and $\l az\r$ and hence is centralised by both $a$ and $az$. Thus $\l a\r\cap\l az\r$ is centralised by $\l a,az\r=\l a,z\r=G$, and so in particular is a normal subgroup of $G$.

The intersection $\l a\r\cap\l zz^a\r$ is contained in the cyclic group $\l a\r$ and so is centralised by $a$. Also the cyclic group $\l zz^a\r$ is inverted by the dihedral group $\l z,z^a\r$, and hence each of its subgroups is normalised by $z$. In particular $\l a\r\cap\l zz^a\r$ is normalised by $z$. Thus $\l a\r\cap\l zz^a\r$ is normalised by $\l a,z\r=G$.

(c)
Let $K_1:= \l az\r$ and $K_2:= \l z,z^a\r$, and for each $i$ set $L_i:= \l a\r\cap K_i$.
Assume that, for some $i$, $|K_i:L_i|\leq 2$. By Lemma~\ref{lem:rotpair}(b), $z, z^a\not\in\l a\r$, and it follows that $L_i\ne K_i$ so  $|K_i:L_i|= 2$, and $z, z^a, az \not\in L_i$.

Suppose first that $i=1$. Then $L_1=\l (az)^2\r$, and by part (b), $L_1\lhd G$. Denote each coset $L_1x$ in $G$ by $\ov x$. Since $G=\l a,z\r$, we have  $G/L_1=\l \ov a,\ov z\r=\l \ov{az},\ov z\r$. Now $\ov{az}$ and $\ov z$ are distinct elements (since $a\not\in L_1$) of order $2$, so $G/L_1$ is a dihedral group $\D_{2r}$, for some $r\geq 2$, with an index $2$ cyclic normal subgroup generated by $\ov{az} \cdot \ov z =\ov a$. Since $L_i\leq \l a\r$ it follows that $\l \ov a\r = \l a\r/L_1$ and hence that $\l a\r\lhd G$.
Thus $G=\l a,z\r=\l a\r{:}\l z\r$, and $\Ga$ has $|G|/|\l a\r|=2$ vertices.
So $\Ga=\K_2^{(\lambda)}$, where $\lambda=|a|$.

Now suppose that $i=2$, so $L_2:=\l a\r\cap K_2$ has index $2$ in $K_2=\l z,z^a\r$, and we note that $L_2$ is cyclic.
If $G$ is abelian then the  conclusion follows from Corollary~\ref{cor:rotpair}(a), so we may assume that $G$ is not abelian, and in particular $z^a\ne z$.
Therefore $K_2 \cong \D_{2r}$ where $r=|zz^a|\geq 2$. If $r=2$ then, since $z, z^a\not\in L_2$, the only possibility for the index $2$ subgroup $L_2$ of  $K_2$ is $\l zz^a\r$, and if $r\geq 3$, then the only cyclic index $2$ subgroup $L_2$ of $K_2$
is  $\l zz^a\r$. Thus in both cases we have $L_2=\l zz^a\r$, and so $L_2=\l a\r\cap\l zz^a\r\lhd G$ by part (b). Denoting each coset $L_2x$ in $G$ by $\ov x$, we have
$G/L_2=\l \ov a,\ov z\r$. Now $\ov z\cdot\ov z^{\ov a}=\ov{zz^a}=1$, and hence
$\ov z^{\ov a}=\ov z$,  so $G/L_2$ is abelian.
In particular, $\l\ov a\r$ is normal in $\ov G$, and since
$\l zz^a\r=L_2\leqslant\l a\r$, it follows that $\l \ov a\r = \l a\r/L_2$ and hence that
$\l a\r$ is normal in $G=\l a,z\r$.
Thus again we find that $G=\l a\r{:}\l z\r$ and  $\Ga$ has $|G|/|\l a\r|=2$ vertices, so $\Ga=\K_2^{(\lambda)}$, where $\lambda=|a|$.
\qed

Using this information we exhibit several cycles in this vertex rotary graph $\Ga$ containing the edge $\l z\r=[\a,e,\b]$. By the \emph{stabiliser in $G$ of a cycle $C$} of $\Ga$, we mean the largest subgroup of $G$ which leaves invariant the sequence class of $C$. The most delicate part of our analysis is proving that the edge sequences we describe involve pairwise distinct edges.  We do not treat the cases where $\Ga$ has only two vertices, or the case where $\Gamma$ is a simple cycle $\C_m$. Equivalent information could be deduced for these graphs using information from Corollary~\ref{cor:rotpair}(b), (c) and Remark~\ref{rem:k2lambda}.

\begin{lemma}\label{lem:cyc-vertrot}
Using the notation from Lemma~$\ref{kernel-V-F}$, suppose that $|V|>2$, equivalently, $\Ga\ne \K_2^{(\lambda)}$, and assume also that $\Ga$ is not a simple cycle.   Then the following hold.
\begin{enumerate}
\item[(a)] Assume that  $m:=|az|$ is finite and let $\lambda':= |\l a\r\cap\l az\r|$. Then $\lambda'$ divides $\gcd(m,\lambda)$ and $m/\lambda'\geq3$. For $i=0,1,\dots,m-1$, the edge
\begin{equation*}
\mbox{$e_i=\l z\r (az)^i$ is of the form $[\l a\r(az)^i,e_i, \l a\r(az)^{i+1}]$,}
\end{equation*}
where $e_0=e$, and the sequence $C(az)=(e_0,e_1,\dots,e_{m-1})$ is an $m$-cycle such that the stabiliser of $C(az)$ in $G$ is $G_{C(az)}=\l az\r\cong Z_m$ and acts regularly and faithfully on the edges of $C(az)$, with $az:e_i\to e_{i+1}$. Moreover  $C(az)$ is a regular cycle with induced subgraph $[C(az)]=\C_{m/\lambda'}^{(\lambda')}$.

\item[(b)] Assume that $\ell:=|zz^a|$ is finite and let  $\lambda'':= |\l a\r\cap\l zz^a\r|$. Then $\lambda''$ divides $\gcd(\ell,\lambda)$ and $\ell/\lambda''\geq3$. For $i=0,1,\dots,\ell-1$, the edge
\begin{align*}
&\mbox{$e_{2i}'=\l z\r (zz^a)^i$ is of the form $[\l a\r z(zz^a)^i,e_{2i}', \l a\r (zz^a)^{i}]$; \ and}\\
&\mbox{$e_{2i+1}'=\l z\r a(zz^a)^i$ is of the form $[\l a\r(zz^a)^i,e_{2i+1}', \l a\r z(zz^a)^{i+1}]$,}
\end{align*}
where $e_0'=e$, and the sequence $C(zz^a)=(e_0',e_1',\dots,e_{2\ell-1}')$  is a $2\ell$-cycle such that the stabiliser of $C(zz^a)$ in $G$ is $G_{C(zz^a)}=\l z,z^a\r\cong \D_{2\ell}$ and acts faithfully on the edges of $C(zz^a)$, with $zz^a:e_{j}'\to e_{j+2}'$, and hence with orbits $\{ e_{2i}'\mid 0\leq i\leq \ell-1\}$ and $\{ e_{2i+1}'\mid 0\leq i\leq \ell-1\}$.  Moreover  $C(zz^a)$ is a regular cycle with induced subgraph $[C(zz^a)]=\C_{2\ell/\lambda''}^{(\lambda'')}$.

%
\end{enumerate}
\end{lemma}

\begin{remark}\label{rem:cyc-vertrot}
{\rm
(a) First we note that Lemma~\ref{lem:cyc-vertrot}(b) does not hold if $\Gamma$ is a simple cycle $\C_n$, and the reasons are different for odd and even values of $n$. If $n$ is odd, then the edge sequence $C(zz^a)$ is not a cycle as each edge is repeated twice, while if $n$ is even, then the subgraph $[C(zz^a)]=\C_{2\ell/\lambda''}^{(\lambda'')}$ is equal to $\Gamma$ (with $\lambda''=1$ and $n=2\ell$) and the full group $G$ preserves the sequence class of $C(zz^a)$ and contains $\l z,z^a\r$ as a subgroup of index $2$.

(b) We obtain two additional cycles containing the edge $e_0$, which are analogous to the cycles $C(az)$ and $C(zz^a)$ by replacing $a$  with $a^{-1}$ in our arguments. In line with the notation in Lemma~\ref{lem:cyc-vertrot} we call these cycles $C(a^{-1}z)$ and $C(zz^{a^{-1}})$. They have the same lengths as $C(az)$ and $C(zz^a)$, respectively.

To see this, observe that $a^{-1}z = (za)^{-1} =((az)^z)^{-1}$ so $|a^{-1}z|=|az|=m$. This means that taking the image of $C(az)$ under $z$ and then reversing the sequence  yields the edge sequence $C(a^{-1}z)$  obtained  by replacing $az$ with
$a^{-1}z$ in Lemma~\ref{lem:cyc-vertrot} (a). (Note that it still contains the edge $e=e^z$.) Hence $C(a^{-1}z)$ is a regular $m$-cycle with stabiliser  $G_{C(a^{-1}z)}=
\l az\r^z=\l za\r=\l a^{-1}z\r \cong \ZZ_{m}$ and induced subgraph $[C(a^{-1}z)]=\C_{m/\lambda'}^{(\lambda')}$.

For the second claim, observe that $zz^{a^{-1}} =zaza^{-1} = za(zz^a)(za)^{-1}$
so $|zz^{a^{-1}}|=|zz^a|=\ell$. This means that taking the image of $C(zz^a)$ under $(za)^{-1}$ yields  the edge sequence $C(zz^{a^{-1}})$ obtained  by replacing $zz^a$ with $zz^{a^{-1}}$ in Lemma~\ref{lem:cyc-vertrot} (b).  Hence $C(zz^{a^{-1}})$  is a regular $2\ell$-cycle with stabiliser  $G_{C(zz^{a^{-1}})}=\l z, z^{a^{-1}}\r\cong \D_{2\ell}$ and induced subgraph $[C(zz^{a^{-1}})]=\C_{2\ell/\lambda''}^{(\lambda'')}$. Note that $e_{1}'$ in $C(zz^a)$ is mapped by $(za)^{-1}$ to
\[
(e_1')^{(za)^{-1}}=  \l z\r a(za)^{-1}= \l z\r = e,
\]
and hence $C(zz^{a^{-1}})$ contains the edge $e$. In Corollary~\ref{cor:cyc-vertrot} we record the details of these additional cycles.

(c) By Lemma~\ref{lem:cyc-vertrot} and part (b), for the cycle $C=C(az)$ or $C(a^{-1}z)$, the stabiliser $G_C$ is transitive on the vertex set and on the edge set, and is bi-regular on the arc set of $[C]$; while for  the cycle $C=C(zz^a)$ or $C(zz^{a^{-1}})$, the stabiliser $G_C$ is transitive on the vertex set, bi-transitive on the edge set, and bi-regular on the arc set of $[C]$.
}
\end{remark}

\begin{corollary}\label{cor:cyc-vertrot}
Using the assumptions and notation from Lemma~\ref{lem:cyc-vertrot}, the following hold.
\begin{enumerate}
\item[(a)] If $m:=|az|$ is finite then $\lambda':= |\l a\r\cap\l az\r|= |\l a\r\cap\l a^{-1}z\r|$, and for $i=0,1,\dots,m-1$, the edge
\begin{equation*}
\mbox{$e_i=\l z\r (a^{-1}z)^i$ is of the form $[\l a\r(a^{-1}z)^i,e_i, \l a\r(a^{-1}z)^{i+1}]$,}
\end{equation*}
where $e_0=e$, and the sequence $C(a^{-1}z)=(e_0,e_1,\dots,e_{m-1})$ is an $m$-cycle such that the stabiliser of $C(a^{-1}z)$ in $G$ is $G_{C(a^{-1}z)}=\l a^{-1}z\r\cong Z_m$ and acts regularly and faithfully on the edges of $C(a^{-1}z)$, with $a^{-1}z:e_i\to e_{i+1}$. Moreover  $C(a^{-1}z)$ is a regular cycle with induced subgraph $[C(a^{-1}z)]=\C_{m/\lambda'}^{(\lambda')}$, and $C(a^{-1}z)$ and $C(az)$ lie in different sequence classes of cycles.

\item[(b)] If $\ell:=|zz^a|$ is finite then  $\lambda'':= |\l a\r\cap\l zz^a\r|= |\l a\r\cap\l zz^{a^{-1}}\r|$, and for $i=0,1,\dots,\ell-1$, the edge
\begin{align*}
&\mbox{$e_{2i}'=\l z\r (zz^{a^{-1}})^i$ is the edge $[\l a\r z(zz^{a{-1}})^i,e_{2i}', \l a\r (zz^{a^{-1}})^{i}]$; \ and}\\
&\mbox{$e_{2i+1}'=\l z\r a^{-1}(zz^{a^{-1}})^i$ is the edge $[\l a\r(zz^{a^{-1}})^i,e_{2i+1}', \l a\r z(zz^{a^{-1}})^{i+1}]$.}
\end{align*}
where $e_0'=e$, and the sequence $C(zz^{a^{-1}})=(e_0',e_1',\dots,e_{2\ell-1}')$  is a $2\ell$-cycle such that the stabiliser of $C(zz^{a^{-1}})$ in $G$ is $G_{C(zz^{a^{-1}})}
=\l z,z^{a^{-1}}\r\cong \D_{2\ell}$ and acts faithfully on the edges of $C(zz^{a^{-1}})$, with $zz^{a^{-1}}:e_{j}'\to e_{j+2}'$, and hence with orbits $\{ e_{2i}'\mid 0\leq i\leq \ell-1\}$ and $\{ e_{2i+1}'\mid 0\leq i\leq \ell-1\}$.  Moreover  $C(zz^{a^{-1}})$ is a regular cycle with induced subgraph $[C(zz^{a^{-1}})]=\C_{2\ell/\lambda''}^{(\lambda'')}$, and $C(a^{-1}z)$ and $C(az)$ lie in different sequence classes of cycles.
\end{enumerate}
\end{corollary}

\noindent
\emph{Proof of Lemma~\ref{lem:cyc-vertrot}.}\quad
(a) Let $L:=\l a\r\cap\l az\r\cong Z_{\lambda'}$. Since $L\leq \l az\r$, the order
$\lambda'=|L|$ divides $m$, and since $L\leq \l a\r$, $L$ fixes $\alpha$. By Lemma~\ref{kernel-V-F}(b), $L$ is normal in $G$, and hence $L$ fixes every vertex. In particular
$L\leq G_{\a\b}$, and by Lemma~\ref{lem:rotpair}(a), $G_{\a\b}=\l a\r\cap \l a\r^z=\l a^k\r\cong Z_\lambda$, so $|L|=\lambda'$ divides $\lambda$. Thus $\lambda'$ divides $\gcd(m,\lambda)$. Further, since  $\Ga\ne \K_2^{(\lambda)}$, it follows from Lemma~\ref{kernel-V-F}(c) that $m/\lambda'=|\l az\r:L|\geq3$.

It follows from the definition of the edges $e_i$, and from Definition~\ref{def:cycles}, that $C(az)$ is an $m$-cycle if and only if $e_0,\dots, e_{m-1}$ are pairwise distinct. Suppose to the contrary that $0\leq i<j\leq m-1$ and $e_i=e_j$, that is to say, $\l z\r (az)^i=\l z\r (az)^j$, or equivalently, $e_0=\l z\r = \l z\r(az)^{j-i}=e_{j-i}$. This implies that $(az)^{j-i}\in\{1,z\}$, and since $0<j-i<m=|az|$, it follows that $(az)^{j-i}=z$ has order $2$. Hence $m$ is even and $j-i=m/2$. Also the pairs of vertices adjacent to these edges are the same, that is, $\{ \l a\r,\l a\r az\}=\{ \l a\r (az)^{m/2},\l a\r (az)^{m/2+1}\}$. Now $z=(az)^{m/2}\not\in\l a\r$ by Lemma~\ref{lem:rotpair}(b), and hence $ \l a\r = \l a\r (az)^{m/2+1}$ and $\l a\r az=\l a\r (az)^{m/2}$. Thus $\l a\r$ contains both $(az)^{m/2+1}$ and $(az)^{m/2-1}$, and hence also $(az)^2=a\cdot a^z$. It follows that $a^z\in\l a\r$ and hence $G_\alpha=\l a\r=\l a\r^z=G_\beta$. By Lemma~\ref{lem:rotpair}(a), $k=1$ which implies that $\Ga= \K_2^{(\lambda)}$, a contradiction. Thus $C(az)$ is an $m$-cycle.

Since $G$ is regular on the arc set of $\Ga$ (by Lemma~\ref{lem:rotpair}(b)) it follows that $G_{C(az)}$ acts faithfully on the edges of $C(az)$. Now $G_{C(az)}$ contains $\l az\r\cong Z_m$ which acts regularly on the edges of $C(az)$, and if $G_{C(az)}$ were strictly larger than
$\l az\r$, then it would contain an element which reverses the arc $(\a,\l z\r,\b)$. The only element that does this is $z$, and if $z\in\l az\r$ then also $a\in\l az\r$, and hence $G=\l a,z\r=\l az\r$, and we have a contradiction by Corollary~\ref{cor:rotpair}(a).  Thus we conclude that  $G_{C(az)}=\l az\r$, and note that $az:e_i\to e_{i+1}$ for all $i$. This implies that $L=G_{C(az),\a}$. Then, writing $\a_i:= \l a\r(az)^i$ for each $i$, the distinct vertices incident with the edges of $C(az)$ are $\a_0, \a_1, \dots, \a_{m/\lambda' -1}$, and we have $\a_i=\a_{i+m/\lambda'}$ for each $i$. It follows that  $[C(az)]=\C_{m/\lambda'}^{(\lambda')}$, and $C(az)$ is a regular cycle.

(b) Proof of the first assertions is similar to the proof in case (a): let $L:=\l a\r\cap\l zz^a\r\cong Z_{\lambda''}$. Since $L\leq \l zz^a\r$, the order
$\lambda''$ divides $\ell$, and since $L\leq \l a\r$, $L$ fixes $\alpha$. By Lemma~\ref{kernel-V-F}(b), $L$ is normal in $G$, and hence $L$ fixes every vertex. In particular
$L\leq G_{\a\b}$, and by Lemma~\ref{lem:rotpair}(a), $G_{\a\b}=\l a\r\cap \l a\r^z=\l a^k\r\cong Z_\lambda$, so $|L|=\lambda''$ divides $\lambda$. Thus $\lambda''$ divides $\gcd(\ell,\lambda)$. Further, since  $\Ga\ne \K_2^{(\lambda)}$, it follows from Lemma~\ref{kernel-V-F}(c) that $\ell/\lambda''=|\l zz^a\r:L|\geq3$. We note, for use below, that this means that $z\not\in \l zz^a\r$ (since $\l z,z^a\r\cong \D_{2\ell}$).

It follows from the definition of the edges $e_i'$, and from Definition~\ref{def:cycles}, that $C(zz^a)$ is a $2\ell$-cycle if and only if $e_0',\dots, e_{2\ell-1}'$ are pairwise distinct. Proving this is somewhat delicate. Suppose to the contrary that $e_u'=e_v'$ with $0\leq u, v\leq 2\ell-1$ and $u\ne v$.
If $u, v$ are both even then we have $(zz^a)^i\in\l z\r$ for some $i$ such that $0<i\leq \ell-1$. Since $|zz^a|=\ell$, this implies that $(zz^a)^i=z$, which as we observed above is not possible. Next assume that $u, v$ are both odd. Then  for some $i$ such that $0<i\leq \ell-1$ we have $a(zz^a)^i\in \l z\r a$. Since $(zz^a)^i\ne 1$, this implies that $a(zz^a)^i=za$. If $i=1$ this equation is equivalent to $aza^{-1}=1$, which is a contradiction. Thus $i\geq 2$ and $(zz^a)^i=z^a$, and hence $e_{2i}=\l z\r (zz^a)^i=\l z\r z^a = \l z\r (zz^a)=e_2$; however we have just shown that this is not possible with $i\geq2$. Thus we must have, say, $u=2i$ and $v=2j+1$ for some $i,j$ between $0$ and $\ell-1$. This implies that $\l z\r (zz^a)^i=\l z\r a(zz^a)^j$ so $(zz^a)^{i-j}\in\{ z, za\}$. In particular $i\ne j$. In the dihedral group
$\l z,z^a\r$, the element $z$ does not lie in $\l zz^a\r$, and hence $(zz^a)^{i-j}=za$, which is equivalent to $(zz^a)^{i-j-1}=az$. This implies that $a\in \l z,z^a\r$. Since $z\not\in\l a\r$ (by Lemma~\ref{lem:rotpair}) it follows that $1\ne az=(zz^a)^{i-j-1}$, and hence $a=(zz^a)^{i-j-1}z$ is an involution in the dihedral group $\l z,z^a\r$ not lying in the cyclic subgroup $\l zz^a\r$. Since $\l a\r$ has order $k\lambda$ this implies that $k\lambda=2$ and we have a contradiction by Corollary~\ref{cor:rotpair}(c), since we are assuming that $\Gamma$ is not $\K_2^{(\lambda)}$ and is not a simple cycle. Thus we have proved that $C(zz^a)$ is a $2\ell$-cycle.

It is clear that the action $zz^a:e_{j}'\to e_{j+2}'$ preserves $C(zz^a)$ and also $z$ induces a reflection of $C(zz^a)$, so the stabiliser $G_{C(zz^a)}$ contains the subgroup $\l z,z^a\r\cong \D_{2\ell}$ acting faithfully on the edges of $C(zz^a)$ with the orbits as in (b). A similar argument to that in part (a) shows that  $C(zz^a)$ is a regular cycle with induced subgraph $[C(zz^a)]=\C_{2\ell/\lambda''}^{(\lambda'')}$.

Finally we prove that $G_{C(zz^a)}=\l z,z^a\r\cong \D_{2\ell}$. If this is not the case then $G_{C(zz^a)}$ contains an element $b$ which induces a rotation of the cycle such that $b:e_i\to e_{i+1}$ for all $i$. In particular $b$ sends $e_0=\l z\r$ to $e_1=\l z\r a$. Hence $ba^{-1} \in G_e=\l z\r$, so $b=z$ or $b=za$. In either case it would follow that  $a\in G_{C(zz^a)}$, since  $b, z\in G_{C(zz^a)}$. Thus $G=\l a,z\r= G_{C(zz^a)}$, and this implies that $G$ is the dihedral group $\D_{4\ell}$ preserving the sequence class of  $C(zz^a)$.  Since $G$ is edge transitive we have $\Gamma=[C(zz^a)]=\C_{2\ell/\lambda''}^{(\lambda'')}$.  In particular the valency $k=2$ and edge multiplicity $\lambda=\lambda''$; and we have shown already that $\ell/\lambda''\geq3$ so $\ell\geq 3\lambda$.
If $a$ lies in the unique cyclic subgroup $Z_{2\ell}$ of $G$ then $z$ inverts $a$ and hence $G=\l a,z\r$ has order $2|a|=2k\lambda=4\lambda$. This is a contradiction since $|G|=4\ell \geq 4\cdot 3\lambda$. Hence $a$ does not lie in this cyclic subgroup and so $|a|=2$, that is $2=k\lambda$, so $\lambda=1$ and $\Gamma$ is a simple cycle, a contradiction. This completes the proof of part (b).
\qed

\noindent
\emph{Proof of Corollary~\ref{cor:cyc-vertrot}.}\quad
By Remark~\ref{rem:core}(b), the only remaining facts to establish are that the cycles $C(a^{-1}z)$ and $C(az)$ lie in different sequence classes, and also $C(zz^{a^{-1}})$ and $C(zz^a)$ lie in different sequence classes.

Suppose that $C(a^{-1}z)$ and $C(az)$ lie in the same sequence class. Now the two edges of $C(a^{-1}z)$ incident with $e=\l z\r$ are $\l z\r a^{-1}z$ and $\l z\r(a^{-1}z)^{-1}=\l z\r a$. Hence $\l z\r az=\l z\r a^{-1}z$ or $\l z\r az=\l z\r a$. The former implies that $a^2\in\l z\r$ which implies that either $a^2=1$ or $z\in\l a\r$, neither of which is possible. The latter implies that $z\in\l z^a\r$ and hence that $z=z^a$. This means that $G$ is abelian and hence  $\Ga=\K_2^{(\lambda)}$ by Corollary~\ref{cor:rotpair}(a), which is a contradiction.

Finally suppose that  $C(zz^{a^{-1}})$ and $C(zz^a)$ lie in the same sequence class. The two edges of $C(zz^a)$ incident with $e=\l z\r$ are $\l z\r a$ and $\l z\r a(zz^a)^{-1}=\l z\r az^az=\l z\r az$. Hence $\l z\r a^{-1}=\l z\r a$ or $\l z\r a^{-1}=\l z\r az$.  The former implies that $a^2\in\l z\r$ which yields a contradiction as in the previous paragraph. Hence the latter holds and so $a^{-1}\in \l z\r az$, so either $a^{-1}=az$ or $a^{-1}=zaz =a^z$. If $a^{-1}=az$ then $z=a^{-2}\in\l a\r$, which is a contradiction. Hence $a^z=a^{-1}$, and this implies that $H^z=H$ and hence that $\Ga$ has valency $k=1$ (by Lemma~\ref{lem:rotpair}) and  so $\Ga=\K_2^{(\lambda)}$, which is a contradiction.
\qed

\section{General constructions of vertex-rotary maps and flag-regular maps}\label{sec:cons}

We first use the cycles identified in Lemma~\ref{lem:cyc-vertrot} and Corollary~\ref{cor:cyc-vertrot} to give two general constructions of vertex-rotary maps. Then, in Subsection~\ref{s:con-regularmap}, we use the theory of coset graphs from Section~\ref{def-coset-g} for a general construction of flag-regular maps.

\subsection{Construction of rotary maps}\label{s:con-rotarymap}

The first construction produces rotary maps using the cycles from Lemma~\ref{lem:cyc-vertrot}(a) and Corollary~\ref{cor:cyc-vertrot}(a).

\begin{construction}\label{cons-rotary}
{\rm
Let $G$ be a group with a rotary pair $(a,z)$ such that $H:=\l a\r$ has order $|a|\geq3$, $H\ne H^z$,  $J:=\l z\r\cong \ZZ_2$, and $|az|$ is finite. Define $\Cos(G,\l a\r,\l z\r,\l az\r)$  as the incidence  configuration  $(V,E,F,\bfI)$,  where
\[
\mbox{$V=[G:\l a\r]$,\quad $E=[G:\l z\r]$, \quad $F=[G:\l az\r]$,}
\]
and the incidence relation $\bfI$ is given by non-empty intersection  between members of $V$, $E$ and $F$. Let
\[
k=|H:H\cap H^z|,\quad \lambda=|H\cap H^z|,\quad m=|az|,\quad \lambda'=|\l a\r\cap\l az\r|,
\]
and associate with $f:=\l az\r\in F$ the edge sequence $C(f):=C(az)=(e_0, e_1,\dots,e_{m-1})$ of Lemma~\ref{lem:cyc-vertrot}(a).  We shall prove that $\Cos(G,\l a\r,\l z\r,\l az\r)$ gives rise to a $G$-rotary map $\RM(G,a,z)$ such that $C(f)$ is the boundary cycle of the face $f$.
}
\end{construction}

The group $G$ acts naturally on the sets $V$, $E$, and $F$  in Construction~\ref{cons-rotary} by right multiplication, and the action preserves the incidence relation $\bfI$, so $G$ is a group of automorphisms of the incidence configuration $\Cos(G,\l a\r,\l z\r,\l az\r)$.  The main result is the following and includes in particular an argument showing how the rotary map $\RM(G,a,z)$ is obtained.

\begin{proposition}\label{p:rotary}
Using the notation from Construction~$\ref{cons-rotary}$,
\begin{enumerate}
\item[(a)] $\Cos(G,\l a\r,\l z\r,\l az\r)$ gives rise to a $G$-rotary map $\RM(G,a,z)$ with valency $k$, edge-multiplicity $\lambda$ and face-length $m$,

\item[(b)] such that, for each $g\in G$, the image of $C(f)$ under $g$ is the boundary cycle of the face $f^g$, and $C(f)$ is a regular $m$-cycle with $[C(f)] =\C_{m/\lambda'}^{(\lambda')}$, where $m/\lambda'\geq3$ and $\lambda'$ divides $\gcd(m,\lambda)$, and $G_f=G_{C(f)}=\l az\r$;

\item[(c)] The kernels of the $G$-actions on $V$ and $V\cup F$ are $G_{(V)}=\l a\r\cap \l a^z\r$ and $G_{(V\cup F)}=\l a\r\cap \l az\r\cong \ZZ_{\lambda'}$. In particular, $\RM(G,a,z)$ is a circular embedding of $\Gamma$ if and only if $G$ acts faithfully on $V\cup F$, or equivalently, if and only if $\lambda'=1$;

\item[(d)]
$z$ interchanges the two faces $f,f^z$ incident with the edge $e=\l z\r$.

\end{enumerate}
\end{proposition}

\proof
By Lemma~\ref{lem:rotpair2}, the underlying graph $\Gamma=(V,E,\bfI)$ is the $G$-vertex-rotary graph $\Cos(G,H,J)$ with valency $k$, edge-multiplicity $\lambda$, and $\Gamma$ is neither $\K_2^{(\lambda)}$ nor a simple cycle, so the results of Section~\ref{s:vertrot} apply. To see that $\Cos(G,\l a\r,\l z\r,\l az\r)$ gives rise to a map with face set $F=[G:\l az\r]$, we first show that each right coset of $\l az\r$ corresponds to a cycle of $\Gamma$. Since $G$ is transitive on $F$ it is sufficient to consider  $f=\l az\r$, and we note that the stabiliser $G_f=\l az\r$ in this action. By definition, it is easily shown that the edges in $E=[G:\l z\r]$ incident with $\l az\r$ are precisely those of the form $\l z\r(az)^i$, for $0\leqslant i\leqslant m-1$, and by Lemma~\ref{lem:cyc-vertrot}(a), these edges may be sequenced to form the regular $m$-cycle $C(f)=C(az)$ of $\Gamma$ with $[\C(f)] =\C_{m/\lambda'}^{(\lambda')}$, where $m/\lambda'\geq3$ and $\lambda'$ divides $\gcd(m,\lambda)$, and moreover the stabiliser in $G$ of the cycle $C(f)$ is equal to $\l az\r$, so $G_f=G_{C(f)}=\l az\r$.

Since $\l a\r$, $\l z\r$ and $\l az\r$ have pairwise non-empty intersections, $(\l a\r,\l z\r,\l az\r)$ is an incident triple. In particular the  edge $e:=\l z\r$ is incident with  $f$ and belongs to $C(f)$ (with $e=e_0$ in Lemma~\ref{lem:cyc-vertrot}(a)). Now $e$ is fixed by $z$ and is incident with $f^z = \l az\r z$; moreover by Remark~\ref{rem:cyc-vertrot}(b), $C(f^z)$, which is the image of $C(az)$ under $z$, lies in the same sequence class as the cycle $C(a^{-1}z)$ of Corollary~\ref{cor:cyc-vertrot}(a). Also, by Corollary~\ref{cor:cyc-vertrot}(a), $C(a^{-1}z)$ and $C(az)$ are not sequence equivalent so the faces $f$ and $f^z$ are distinct. We claim that $f$ and $f^z$ are the only cosets in $F$ incident with $e$. For suppose that $f^g=\l az\r g$ is incident with $e$, for some $g\in G$.
Then $\l z\r\cap\l az\r g\not=\emptyset$ and so, for some integer $i$, $(az)^ig=1$
or $z$. In the former case, $g=(az)^{-i}$, and hence $f^g=\l az\r g=\l az\r = f$, while in the latter case, $g=(az)^{-i}z$ and $f^g=\l az\r g = \l az\r z=f^z$. This proves the claim.
Thus each edge of $\Ga$ is incident with exactly two cosets in $F$, and hence lies in exactly two of the corresponding boundary cycles.
In particular, $z$ interchanges the two faces $f,f^z$ incident with $e=\l z\r$, proving part~(d).

To complete the proofs of parts (a) and (b), we identify each cycle $C(f^g)$, for $g\in G$, with the boundary of a Euclidean disc $\D(f^g)$, and
then `sew' these discs together by identifying arcs on two such discs when they correspond to the same arc of the corresponding cycles. This defines a topological space  $\calS$ with the graph $\Ga$ embedded in it.
To verify that $\calS$ is a surface, we prove:

\smallskip
\noindent
{\bf Claim:}\ {\it each point $\o$ in $\calS$ has an open neighbourhood homeomorphic to an open disc.}

\noindent
\medskip
The claim is clearly true for each point of a disc $\D(f^g)$ which is not on the boundary $C(f^g)$  (for some $g\in G$).
Since an edge in $E$ is incident with exactly two such cycles, the claim is also true for each interior point of an edge. It remains to consider the case where $\o$ is a vertex in $V$. Since the action of $G$ on $V, E, F$ is induced by some group of homeomorphisms of $\calS$, and since $G$ is transitive on $V$, it is sufficient to consider $\o=\a$.
The regular action of  $G_\a=\l a\r$ on $E(\a)$ arranges the edges in $E(\a)$ in order as
$e,e^a,e^{a^2},\dots,e^{a^{k\lambda-1}}$.
Let $C_0=C(f)$, and $C_i=C_0^{a^i}$, where $0\leqslant i\leqslant k\lambda-1$, and for each $i$ let $\D_i$ be the disc of $\calS$ corresponding to $C_i$.
Then for each $i$, $C_i$ contains the 2-arc $(e^{a^{i}},\a,e^{a^{i+1}})$.
These $k\lambda$ discs around $\a$ form a larger disc $\D_0\cup\dots\cup\D_{k\lambda-1}$ with $\a$ being an interior point. This proves the claim.
Hence $\calS$ is a surface, giving rise to the map $(V,E,F)=\RM(G,a,z)$, which is a $G$-rotary map by  definition. Thus parts (a) and (b) are proved.

By Lemma~\ref{kernel-V-F}(a),  $G_{(V)}= \l a\r\cap \l a^z\r$.  The kernel $G_{(V\cup F)}$ fixes in particular the vertex $\a=\l a\r$ and the face $\l az\r$,
and hence $G_{(V\cup F)}\leq \l a\r\cap \l az\r$. By Lemma~\ref{kernel-V-F}(b), $\l a\r\cap\l az\r$ is normal in $G$, and hence $\l a\r\cap\l az\r$ fixes every vertex and every face. So equality holds and we have $G_{(V\cup F)}= \l a\r\cap \l az\r\cong \ZZ_{\lambda'}$. The  boundary cycles $C(f)$ are simple cycles if and only if $\lambda'=1$, and thus by Definition~\ref{def:circembedding}, $\RM(G,a,z)$ is a circular embedding  if and only if $\lambda'=1$, or equivalently, if and only if $G$ is faithful on $V\cup F$.
\qed

\subsection{Construction of bi-rotary maps}\label{s:con-bi-rotarymap}

The second construction uses the cycles identified in Lemma~\ref{lem:cyc-vertrot}(b) and Corollary~\ref{cor:cyc-vertrot}(b) and produces bi-rotary maps.
In this case we cannot simply give the construction as a coset configuration (see Remark~\ref{rem:bi-rotary}). Instead we identify the faces with their boundary cycles and define an edge and face to be incident if the edge belongs to the corresponding cycle, and similarly a vertex and face are incident if the vertex is incident with some edge in the corresponding cycle.

\begin{construction}\label{cons-bi-rotary}
{\rm
Let $G$ be a group with a rotary pair $(a,z)$ such that $H:=\l a\r$ has order $|a|\geq3$, $H\ne H^z$,  $J:=\l z\r\cong \ZZ_2$, $|zz^a|$ is finite, and let $W:=\l z,z^a\r$.  Associate with $W$ the edge sequence $C(W):=C(zz^a)=(e_0', e_1',\dots,e_{2\ell-1}')$ of Lemma~\ref{lem:cyc-vertrot}(b). Define $\Cos(G,H, J, W)$  as the incidence  configuration  $(V,E,F,\bfI)$,  where
\[
\mbox{$V=[G:\l a\r]$,\quad $E=[G:\l z\r]$, \quad $F=\{ C(W)g\mid g\in G\}$,}
\]
and the incidence relation $\bfI$ is given by non-empty intersection  between members of $V$ and $E$, inclusion between members of $E$ and $F$,  and incidence as described above between members of $V$ and $F$. Let
\[
k=|H:H\cap H^z|,\quad \lambda=|H\cap H^z|,\quad \ell=|zz^a|,\quad \lambda''=|\l a\r\cap\l zz^a\r|.
\]
We shall prove that $\Cos(G,H, J, W)$ gives rise to a $G$-bi-rotary map $\BRM(G,a,z)$.
}
\end{construction}

The group $G$ acts naturally on the sets $V$, $E$, and $F$  in Construction~\ref{cons-bi-rotary} by right multiplication, and the action preserves the incidence relation $\bfI$, so $G$ is a group of automorphisms of the incidence configuration $\Cos(G,H,J,W)$.  The main result is the following and includes  an argument similar to that in Proposition~\ref{p:rotary} to show how the bi-rotary map $\BRM(G,a,z)$ is obtained.

\begin{proposition}\label{p:bi-rotary}
Using the notation from Construction~$\ref{cons-bi-rotary}$,
\begin{enumerate}
\item[(a)] $\Cos(G,H,J,W)$ gives rise to a $G$-bi-rotary map $\BRM(G,a,z)$ with valency $k$, edge-multiplicity $\lambda$ and face-length $2\ell$,

\item[(b)] such that, for each $g\in G$, the image of $C(W)$ under $g$ is (the boundary cycle of) the face $C(W)g$, and $C(W)$ is a regular $2\ell$-cycle with $[\C(W)] =\C_{2\ell/\lambda''}^{(\lambda'')}$, where $\ell/\lambda''\geq3$ and $\lambda''$ divides $\gcd(\ell,\lambda)$, and $G_{C(W)}=\l z,z^a\r = W$.

\item[(c)] The kernels of the $G$-actions on $V$ and $V\cup F$ are $G_{(V)}=\l a\r\cap \l a^z\r$ and $G_{(V\cup F)}=\l a\r\cap \l zz^a\r$. In particular $\BRM(G,a,z)$ is a circular embedding of $\Gamma$ if and only if $G$ acts faithfully on $V\cup F$, if and only if $\lambda''=1$.

\item[(d)]
$z$ fixes each of the two faces incident with the edge $e=\l z\r$.

\end{enumerate}
\end{proposition}

\proof
By Lemma~\ref{lem:rotpair2}, the underlying graph $\Gamma=(V,E,\bfI)$ is the $G$-vertex-rotary graph $\Cos(G,H,J)$ with valency $k$, edge-multiplicity $\lambda$, and $\Gamma$ is neither $\K_2^{(\lambda)}$ nor a simple cycle, so the results of Section~\ref{s:vertrot} apply. By Lemma~\ref{lem:cyc-vertrot}(b), $C(W)=C(zz^a)$ is a regular $2\ell$-cycle
of $\Gamma$ with $[C(W)] =\C_{2\ell/\lambda''}^{(\lambda'')}$, where $\ell/\lambda''\geq3$ and $\lambda''$ divides $\gcd(\ell,\lambda)$, and moreover the stabiliser in $G$ of the cycle $C(W)$ is equal to $\l z,z^a\r$, so $G_{C(W)}=\l z,z^a\r=W$. This proves part (b).

To see that $\Cos(G,H,J,W)$ gives rise to a map with face set $F$, we next show that
each edge is incident with exactly two cycles in $F$. Since $G$ acts transitively on the edge set $E$, it is sufficient to prove this for the edge $e=\l z\r$.
By the definition of incidence it follows that $(\l a\r,\l z\r,C(W))$ is an incident triple. In fact $e:=\l z\r$ is the  edge $e_0'$ of the cycle $C(W)=C(zz^a)$ of Lemma~\ref{lem:cyc-vertrot}(b). By Remark~\ref{rem:cyc-vertrot}(b), $e$ is also incident with the cycle $C(zz^{a^{-1}})$ of Corollary~\ref{cor:cyc-vertrot}(b), which is the image of $C(W)$ under $(za)^{-1}$, and by Corollary~\ref{cor:cyc-vertrot}(b), these two cycles
$C(W)$ and $C(W)(za)^{-1}$ are distinct elements of $F$. Suppose that $e$ is incident with $C(W)g\in F$, for some $g\in G$. Then there is an integer $i$ such that
$\l z\r = \l z\r (zz^a)^i$ or $ \l z\r a(zz^a)^i$. Thus $(zz^a)^ig\in\l z\r$ or $a(zz^a)^ig\in\l z\r$, respectively. If $(zz^a)^ig\in\l z\r$, then $g\in\l z, zz^a\r\leq
\l z,z^a\r=W=G_{C(W)}$, so $C(W)g=C(W)$. So suppose that $a(zz^a)^ig\in\l z\r$, that is, $a(zz^a)^ig=z^\delta$ for $\delta=0$ or $1$, and so $g=(zz^a)^{-i}a^{-1}z^\delta$. Note that $(zz^a)^{-i}\in W=G_{C(W)}$, and hence $C(W)g=C(W)a^{-1}z^\delta$. If $\delta=1$ then $C(W)g=C(W)a^{-1}z=C(W)(za)^{-1}$. So assume that $\delta=0$. Then $C(W)g=C(W)a^{-1} = C(W) z^a a^{-1}$ (since $z^a\in W=G_{C(W)}$), and this equals $C(W)a^{-1}z = C(W)(za)^{-1}$.
Thus $e$, and hence also each edge of $\Ga$, is incident with exactly two cycles in $F$.
By Lemma~\ref{lem:cyc-vertrot}(b), $z$ fixes $C(zz^a)= C(W)$, and since  $e=\l z\r$ is incident with exactly two faces in $F$, namely $C(W)$ and $C(W)(za)^{-1}$, it follows that $z$ fixes each of these faces.

Note that, as $G_{C(W)}=W$, the cycles $C(W)g$ in $F$ depend only on the right coset $Wg$ containing $g$. To complete the proof of part (a), we identify each cycle $C(W)g\in F$ with the boundary of a Euclidean disc $\D(Wg)$, and
then `sew' these discs together by identifying arcs on two such discs when they correspond to the same arc of the corresponding cycles. This defines a topological space  $\calS$ with the graph $\Ga$ embedded in it. Proof that $\calS$ is a surface is exactly the same as that given for the `Claim' in the proof of Proposition~\ref{p:rotary}, so we omit the details. Thus we obtain the map $(V,E,F)=\BRM(G,a,z)$. This map is $G$-vertex-rotary but not $G$-face-rotary (since $G_{C(W)}=W$ is not cyclic). The involution $z$ interchanges the two arcs $(\l a\r,\l z\r,\l a\r z)$ and $(\l a\r z,\l z\r,\l a\r)$, and fixes each of the two faces $C(W)$ and $C(W)(za)^{-1}$ incident with the edge $\l z\r$.
Thus the local orientations of the surface $\calS$ induced by the cyclic actions of
$\l a\r$ on the faces incident with $\l a\r$, and  of $\l a^z\r$ on the faces incident with $\l a\r z$ are different. Therefore the map  $\BRM(G,a,z)$ is $G$-bi-rotary, as in Definition~\ref{def:rotary}(c).

By Lemma~\ref{kernel-V-F},  $G_{(V)}= \l a\r\cap \l a^z\r$.
The kernel $G_{(V\cup F)}$ fixes in particular the vertex $\a=\l a\r$ and the face $C(W)$,
and since $G_{C(W)}=W$ it follows that $G_{(V\cup F)}\leq \l a\r\cap W$. Hence $G_{(V\cup F)}$ is a cyclic normal subgroup of $W=\l z,z^a\r\cong \D_{2\ell}$. By part (b), $\ell\geq \ell/\lambda''\geq 3$, so each cyclic normal subgroup of $W$ is contained in $\l zz^a\r$. Thus $G_{(V\cup F)}\leq \l a\r\cap \l zz^a\r$. On the other hand,  by Lemma~\ref{kernel-V-F}(b), $\l a\r\cap \l zz^a\r$ is normal in $G$ and hence (since $G_\a = \l a\r$ and $zz^a\in G_{C(W)}=\l z,z^a\r$),   $\l a\r\cap \l zz^a\r$ acts trivially on $V\cup F$. We conclude that $G_{(V\cup F)}= \l a\r\cap \l zz^a\r\cong \ZZ_{\lambda''}$. Finally,  by Definition~\ref{def:circembedding}, $\RM(G,a,z)$ is a circular embedding  if and only if the  cycles $C(W)g$ are simple cycles, and this holds by part (b) if and only if $\lambda''=1$, or equivalently, if and only if  $G_{(V\cup F)}=1$.
\qed

\begin{remark}\label{rem:bi-rotary}
{\rm
We have seen in Proposition~\ref{p:bi-rotary} that the $G$-action on the set of faces of $\BRM(G,a,z)$ is equivalent to its right multiplication action on $[G:W]$ with $W$ the stabiliser of the face $C(W)$. This suggests a natural question concerning Construction~\ref{cons-bi-rotary}: whether the map could be defined in a simpler way as a coset configuration, similar to that in Construction~\ref{cons-rotary}. Namely, could one take the set of faces to be $[G:W]$ with incidence between elements of $V, E, F$ given by nontrivial intersection. It turns out that this is not possible, and we give in Proposition~\ref{p:birotary-ex} an example of two different $G$-bi-rotary maps arising from the same subgroups $H, J, W.$ Although the group $J$ determines the involution $z$, the subgroup triple $H, J, W$ does not determine the rotary pair $(a,z)$ and does not determine the pair of boundary cycles incident with a given edge (not even up to sequence equivalence).
}
\end{remark}

\begin{proposition}\label{p:birotary-ex}
The group $G=3.A_6$ (non-split central extension) has rotary pairs $(a,z)$ and $(a',z)$ such that
\begin{itemize}
\item[(a)] $(\l a\r,\l z\r,\l z,z^a\r)=(\l a'\r,\l z\r,\l z,z^{a'}\r)$, and
\item[(b)] $\BRM(G,a,z)\not=\BRM(G,a',z)$.
\end{itemize}
In particular the claims of Corollary~$\ref{cor:iso}$ are valid.
\end{proposition}

\proof
The group $G$ is perfect with centre $\Z(G)=\ZZ_3$.
Let $\ov G=G/\Z(G)$, so $\ov G=\A_6$ and $\ov G$ acts naturally on $\{1,2,3,4,5,6\}$.
There exist unique elements $b, z\in G$ such that $|b|= 5$, $|z|=2$, and $\ov b=(12345)$, $\ov z=(34)(56)$. Clearly  $\l \ov b,\ov z\r$ is a $2$-transitive subgroup of $\ov G = A_6$, and an easy calculation shows that $\ov{b^2 z} = (14)(2365)$, of order $4$; hence
$\l \ov b,\ov z\r=\ov G=\A_6$. Since $G$ is a non-split extension it follows that $G=\l b,z\r$. Another easy calculation shows that $\ov{zz^b}=(16435)$.

Let $\Z(G)=\l c\r$, and $a=b^2c$. Then $|a|=15$, and $a^3=b$. Also  $\l a,z\r=G$, so $(a,z)$ is a rotary pair for $G$. Now $\ov{zz^a}=\ov{zz^b}=(16435)$ has order 5, and we claim that $|zz^a|=5$. If this is not the case then $|zz^a|=15$ and we have
$(zz^a)^5\in\Z(G)$; however $z$ inverts $zz^a$, and hence $z$ also inverts $(zz^a)^5$, contradicting the fact that $(zz^a)^5\in\Z(G)$. This proves the claim and so
$\l z,z^a\r=\D_{10}$. Thus in the $G$-bi-rotary map $\BRM(G,a,z)$ we have
$G_{(V\cup F)}=\l a\r\cap\l z,z^a\r=1$, and so $\BRM(G,a,z)$ is a circular map, by Proposition~\ref{p:bi-rotary}(c).

Now consider $a'=a^{11}$.
Then $|a'|=|a|=15$, and $\l a\r=\l a^{11}\r$, so $(a',z)$ is also a rotary pair for $G$.
Moreover, $a^{10}=(b^2c)^{10}=c\in\Z(G)$ centralizes $z$.
Thus $z^{a'}=z^{a^{11}}=z^a$, and so $\l z,z^a\r=\l z,z^{a^{11}}\r$, and we have $(\l a\r,\l z\r,\l z,z^a\r)=(\l a'\r,\l z\r,\l z,z^{a'}\r)$, as in part (a).

In $\BRM(G,a,z)$, the two faces incident with the arc $(\l a\r,\l z\r)$ are uniquely determined by the two 2-arcs $(\l z\r a^{-1},\l a\r, \l z\r)$ and $(\l z\r a,\l a\r, \l z\r)$ of the underlying graph $\Gamma=\Cos(G,H,J)$. Similarly,
in the map $\BRM(G,a',z)$, the two faces incident with the arc $(\l a\r,\l z\r)$ are uniquely determined by the two 2-arcs $(\l z\r a^{-11},\l a\r, \l z\r)$ and $(\l z\r a^{11},\l a\r, \l z\r)$ of the same underlying graph $\Gamma$.  Thus if the two maps $\BRM(G,a,z)$ and $\BRM(G,a^{11},z)$ were equal, then $\l z\r a^{11}, \l z\r$ would be consecutive edges in one of the two cycles containing $\l z\r$ in $\BRM(G,a,z)$. We note that these cycles are determined up to sequence equivalence, and that these two cycles are $C(zz^a)$ and $C(zz^{a^{-1}})$ from Lemma~\ref{lem:cyc-vertrot} and Corollary~\ref{cor:cyc-vertrot}. The pairs of edges consecutive with $\l z\r$ in these cycles are
$\{ \l z\r az, \l z\r a\}$ and  $\{\l z\r a^{-1}z, \l z\r a^{-1}\}$  respectively. Thus $\l z\r a^{11}$ is one of these four edges. If $\l z\r a^{11}= \l z\r a^{\delta}$ with $\delta=\pm 1$, then $\l z\r$ would contain $a^{11\pm 1}$, which is a contradiction since $|a^{10}|=3$ and $|a^{12}|=5$. Thus $\l z\r a^{11}= \l z\r a^{\delta}z$ with $\delta=\pm 1$, which implies that $a^{11}za^{\pm 1}\in\l z\r$. Computing in $\ov G$, and noting that $\ov{a^{11}}=\ov a = \ov{b^2}=(13524)$ we find that
\[
\ov{a^{11}za^{-1}} = (12)(36),\quad \mbox{and}\quad \ov{a^{11}za} = (36254),
\]
neither of which lies in $\l \ov z\r$. Thus part (b) is proved and hence also Corollary~$\ref{cor:iso}$.
\qed

\subsection{Construction of flag-regular maps}\label{s:con-regularmap}

The third construction is also a coset configuration, and uses a coset graph for a group with a flag-regular triple which is defined as follows.

\begin{definition}\label{def-flagregtriple}
{\rm
For pairwise distinct involutions $x,y,z$ in a group $G$, the ordered triple $(x,y,z)$ is called a \emph{flag-regular triple} for $G$ if
\[
\mbox{$G=\l x,y,z\r$, $xz=zx$, $z\not\in\l x,y\r$, and the orders $|xy|$ and $|yz|$ are finite.}
\]
}
\end{definition}

\begin{construction}\label{cons-flagregular}
{\rm
Let $G$ be a group with a flag-regular triple $(x,y,z)$, and let $H:=\l x,y\r$, $J:=\l x,z\r\cong\ZZ_2^2$, and $W:=\l y,z\r$. Define $\Cos(G,H,J,W)$ as the incidence
configuration $(V,E, F,\bfI)$, where
\[
\mbox{$V=[G:H]$,\quad $E=[G:J]$, \quad $F=[G:W]$,}
\]
and the incidence relation $\bfI$ is given by non-empty intersection  between members of $V$, $E$ and $F$. Let $a:=xy$ and $b:= zy$, and let $k$ divide $|a|$ such that
\[
 \l a\r\cap \l a^z\r =\l a^k\r,\quad \lambda = |a^k|,\quad m=|b|,\quad \lambda'=|\l a\r\cap\l b\r|.
\]
We shall prove that $\Cos(G,H, J, W)$ gives rise to a locally finite $G$-flag-regular map $\FRM(G,x,y,z)$ whenever $k\lambda\geq3$ and $m\geq3$.
}
\end{construction}

We note that the assumption $z\not\in\l x,y\r$ is equivalent to the condition that $H$ is a proper subgroup of $G$. It also implies that $H\cap J=\l x\r$ and hence that $|J:H\cap J|=2$ so that all the conditions of Construction~\ref{def-m-graph} are satisfied, and the coset graph $\Cos(G,H,J)$ is well defined.  As in our investigation of vertex-rotary graphs in Lemma~\ref{lem:cyc-vertrot}, and hence in our constructions of vertex-rotary maps, we avoid the degenerate cases where this graph is $\K_2^{(\lambda)}$ or a simple cycle; and we note that these cases will be examined in detail in \cite{Cycles2}.

\begin{lemma}\label{lem:regmapgraph}
Let $G, H, J, W, a, b, k, m, \lambda, \lambda'$ be as in Construction~$\ref{cons-flagregular}$, and let $\Gamma$ be the coset graph $\Cos(G, H, J)$. Assume that $\Gamma$ is not $\K_2^{(\lambda)}$ and is not a simple cycle. Then the following statements hold, where
$\alpha=H\in V$, $\beta=Hz\in V$, and $e=J$ is the edge $[\alpha, e,\beta]=[H,J,Hz]$ of $\Gamma$.
\begin{enumerate}
\item[(a)] $\Gamma$ is a connected $G$-arc-transitive graph of valency $k\geq2$ with edge-multiplicity $\lambda$ such that $|E(\alpha)|=k\lambda\geq3$. Also $\Core_G(H\cap J)=1$ and $G\leq \Aut\Gamma$ (with its action on right cosets), and $G_{\alpha \beta}=
H\cap H^z = \l a^k,x\r\cong \D_{2\lambda}$.

\item[(b)]  $m=|b|=|zy|\geq3$ and  $x\not\in\l y,z\r$.

\item[(c)] $\lambda'$ divides $\gcd(m,\lambda)$ and $m/\lambda'\geq2$. For $i=0,1,\dots,m-1$, the edge
\begin{equation*}
\mbox{$e_i=J b^i$ is of the form $[Hb^i,e_i, Hb^{i+1}]$,}
\end{equation*}
where $e_0=e$, and the sequence $C(W)=(e_0,e_1,\dots,e_{m-1})$ is an $m$-cycle such that the stabiliser of $C(W)$ in $G$ is $G_{C(W)}=W\cong \D_{2m}$ and acts faithfully and arc-regularly on $C(W)$, with $b:e_i\to e_{i+1}$ and $z:e_i\to e_{m-i}$. Moreover $C(W)$ is a regular $m$-cycle and the induced subgraph $[C(W)]$ is either $\C_{m/\lambda'}^{(\lambda')}$ with $m/\lambda'\geq3$, or  $\K_{2}^{(m)}$ with $m=2\lambda'\geq4$.

\item[(d)] $N:=\l a\r\cap\l a^z\r=\l a^k\r$ is normal in $G$ and of index $2$ in $H\cap H^z$; moreover either
	\begin{enumerate}
	\item[(i)]  $G_{(V)}=N$, or
	\item[(ii)] $k=2, \lambda>1, m/\lambda'\geq3$, $\Gamma=\C_{m/\lambda'}^{(\lambda)}$, $G_{(V)}=H\cap H^z=N\rtimes \l x\r\cong \D_{2\lambda}$, and $G/G_{(V)}\cong\D_{2m/\lambda'}$.

	\end{enumerate}

\end{enumerate}
\end{lemma}


\proof
(a)  By Theorem~\ref{thm:cosetgraph}, $\Gamma$ is $G$-arc-transitive of valency $k$ and with edge-multiplicity $\lambda$, and $\Gamma$ is connected since $G=\l x,y,z\r$. Also  $|E(\alpha)|=k\lambda$ by Lemma~\ref{coset-graphs}(ii). If $k=1$ then $\Gamma=\K_2^{(\lambda)}$ since $\Gamma$ is connected, and we have excluded this case, so $k\geq2$. If $|E(\alpha)|=k\lambda=2$ then $k=2, \lambda=1$ and $\Gamma$ is a simple cycle, which is also excluded. Hence $k\lambda\geq3$.
Now $H\cap J=\l x\r\cong\ZZ_2$. If $\Core_G(H\cap J)\ne1$, then $\Core_G(H\cap J)=\l x\r$ is normal in $G$, and in fact is central in $G$. This would imply that $x$ is a central involution of the dihedral group $H$, whereas $x$ inverts $a$. Hence $k\lambda=|a|=2$, which is a contradiction.  Thus $\Core_G(H\cap J)=1$, so $G$ acts faithfully on $\Gamma$ by Theorem~\ref{thm:cosetgraph}. Since $\beta=Hz$, it follows that $G_\beta=H^z$ and
$G_{\alpha\beta}=H\cap H^z$ has index $k$ in $H$. Then as $x\in H\cap J$ fixes $\alpha$ and $e=J$, it follows that $x$ also fixes the second vertex $\beta$ incident with $e$, so $x\in G_{\alpha\beta}$. Therefore $G_{\alpha\beta}=\l a^k,x\r\cong \D_{2\lambda}$, and part (a) is proved.

(b) Since $z\ne y$, we have $m=|zy|\ne 1$. Suppose that $m=2$. Then $y$ commutes with $z$, and by Definition~\ref{def-flagregtriple}, also $x$ commutes with $z$, so $z$ is central in $G$. Since $z\not\in\l x,y\r=H$, by Definition~\ref{def-flagregtriple},  it follows that $G=H\times\l z\r$, and hence $|V|=|G:H|=2$. This however implies that $\Gamma=\K_2^{(\lambda)}$  and we have excluded this case. Thus $m\geq3$. Next suppose that $x\in\l y,z\r=W$. Then $G=W$ and $x$ is an involution in the dihedral group $W\cong \D_{2m}$ (with $m\geq3$) such that $x$ commutes with $z$. If $m$ were odd then $C_W(z)=\l z\r$ and since $z\ne x$ we have a contradiction. Thus $m=2s$ with $s\geq2$, and $C_W(z)=\l z, b^s\r\cong\ZZ_2^2$, yielding $x\in \{b^s, zb^s\}$. If $x=b^s$ then $\l x\r=H\cap J$ would be central in $G$ contradicting part (a). Hence $x=zb^s$.  In this case $a=xy = zb^sy=zyb^s=b^{s+1}$ and $\gcd(s+1, 2s)=\gcd(s+1,2)$. Thus $\l a\r=\l b\r$ if $s$ is even, and $\l a\r=\l b^2\r$ if $s$ is odd. In either case, $H=\l x,y\r$ contains $y$ and $b^2$, and $\l y,b^2\r$ has index at most $2$ in $W=\l y,z\r=G$.
Thus $|V|=|G:H|\leq 2$, which is a contradiction. We conclude that $x\not\in\l y,z\r$ and part (b) is proved.

(c)  Let $L:=\l a\r\cap\l b\r$. Then $L\cong Z_{\lambda'}$, by the definition of $\lambda'$.
Since $L\leq \l b\r$, the order
$\lambda'=|L|$ divides $m$, and since $L\leq \l a\r$, $L$ fixes $\alpha$.
Now $L$ is normalised by $x, y$ (since $L\leq \l a\r$) and by $z$ (since  $L\leq \l b\r$), and hence $L$ is normal in $G$. Then, as $L$ fixes $\alpha$ it follows that $L$ fixes $V$ pointwise. In particular $L\leq \l a\r\cap G_{\alpha\beta}=\l a^k\r\cong\ZZ_\lambda$ by part (a), so $\lambda'$ divides $\lambda$. Thus  $\lambda'$ divides $\gcd(m,\lambda)$. Suppose next that $m=\lambda'$. Then $L=\l b\r$, and we have $b\in G_\alpha=H$ as well as $y\in H$ by definition. Hence $\l y,z\r=W=\l b,y\r\leq H$, and this implies that $z\in H$, contradicting Definition~\ref{def-flagregtriple}. Thus $m/\lambda'\geq2$.

It follows from the definition of the edges $e_i$, and from Definition~\ref{def:cycles}, that $C(W)$ is an $m$-cycle if and only if $e_0,\dots, e_{m-1}$ are pairwise distinct. Suppose to the contrary that $0\leq i<j\leq m-1$ and $e_i=e_j$, that is to say, $J b^i=J b^j$, or equivalently, $e_0=J = Jb^{j-i}=e_{j-i}$. This implies that $b^{j-i}\in J=\{1,x,z,xz\}$, and since $0<j-i<m=|b|$, it follows that $b^{j-i}$ has order $2$. Hence $m$ is even, say $m=2s$, and $b^{j-i}=b^s$ is a central involution in $W$. By part (b), as $m\geq3$, we have $\l b^s\r= Z(W)$. If $b^s=x$ then $\l x\r=H\cap J$ would be normal in $G$ contradicting part (a). Hence $b^s\in\{z, xz\}$. If $b^s=z$ then $z\in Z(W)$ contradicting the fact that $z$ inverts $b$ and $m\geq3$. Thus $b^s=xz$. This implies that $xz\in W$ and hence $x\in W$, again contradicting part (b). Thus the edges $e_i$ are pairwise distinct and $C(W)$ is an $m$-cycle. By the definition of $C(W)$ we have
$b:e_i\to e_{i+1}$ for each $i$. Also $z$ maps $e_i=Jb^i$ to $Jb^iz=Jzb^iz=J(b^i)^z=Jb^{-i}= e_{m-i}$, so $z:e_i\to e_{m-i}$. Thus the stabiliser of $C(W)$ in $G$ contains $\l b,z\r=W$ acting faithfully and arc-regularly on $C(W)$ as $\D_{2m}$. If the stabiliser were strictly larger than $W$ then it would contain a nontrivial element fixing the arc $(\alpha, e_0,\beta)$. However the stabiliser of this arc is $H\cap J=\l x\r$, and this would imply that $x$ stabilises $C(W)$, and hence that $G=\l x,y,z\r$ stabilises $C(W)$. By part (a), $G$ acts faithfully on $V\cup E$ and  this implies that $x\in G\leq \l b,z\r$ which is a contradiction. Thus $G_{C(W)}=W$.

Finally we consider the edge-induced subgraph $[C(W)]$. As $\l b\r\cong\ZZ_m$ acts regularly on the edges of this subgraph, $[C(W)]$ has $m/\lambda'=|\l b\r:\l a\r\cap \l b\r|$  vertices. If $m/\lambda'\geq3$, then it follows from \cite[Theorem 1.1(I)]{Cycles} that $[C(W)]=\C_{m/\lambda'}^{(\lambda')}$, while if $m/\lambda'=2$ then $[C(W)]=\K_{2}^{(m)}$ and as $m=2\lambda'$ is even and $m\geq3$ by part (b), we have $m\geq 4$ (and this subgraph does indeed admit a symmetrical Euler cycle, see
 \cite[Proposition 2.2(b)]{Cycles}). In either of these two cases, $C(W)$ is a regular cycle, by Definition~\ref{def:cycles}.

 (d) Let $N:=\l a\r\cap\l a^z\r=\l a^k\r
\cong \ZZ_\lambda$, the cyclic index $2$-subgroup of $G_{\alpha\beta}=H\cap H^z$. Now $N$ is normalised by $x$ and $y$ (since $N\leq \l a\r$) and also by $z$ (since $z$ interchanges $\l a\r$ and$\l a^z\r$). Hence $N$ is normal in $G$ and so are all subgroups of $N$. Since $N$ fixes the vertex $\alpha$ it follows that $N$ acts trivially on $V$, that is, $N\leq G_{(V)}$. On the other hand, $G_{(V)}\leq G_{\alpha\beta}=H\cap H^z = N\rtimes \l x\r$.

Suppose that $N\ne G_{(V)}$. Then $G_{(V)}=N\rtimes \l x\r\cong \D_{2\lambda}$. If $N=1$ this implies that $\l x\r\lhd\ G$ and hence $\Core_G(H\cap J)=\l x\r$, contradicting Lemma~\ref{lem:regmapgraph}(a). Hence $N\ne 1$, that is, $\lambda>1$. Then $G_{(V)}=N\rtimes \l x\r\cong \D_{2\lambda}$ is normal in $H\cong \D_{2k\lambda}$ and it follows that $k=2$ and, since $\Gamma\ne\K_2^{(\lambda)}$, the base graph of $\Gamma$ is a cycle $\C_n$ and $\Gamma=\C_n^{(\lambda)}$, for some $n\geq3$ by Theorem~\ref{thm:cosetgraph}(d). Noting that $y, z$ act nontrivially on $V$, $\overline{G}:=G/G_{(V)}$ is generated by involutions $\overline{y}=G_{(V)}y$ and  $\overline{z}=G_{(V)}z$, and is isomorphic to $\Aut\C_n = \D_{2n}$. On the other hand $\l b\r=\l zy\r$ acts on $\C_n$ as the rotation group of order $n$, and by part (c) it follows that $\l b\r\cap G_{(V)}=L\cong\ZZ_{\lambda'}$, so $m/\lambda'=n\geq 3$, and $\overline{G}\cong \D_{2m/\lambda'}$.
\qed

Now we verify that Construction~\ref{cons-flagregular} yields a locally finite flag-regular map. We exclude the cases where $\Gamma$ is  $\K_2^{(\lambda)}$ or $\C_n^{(\lambda)}$ for some $n\geq3$.

\begin{proposition}\label{p:flagregular}
Using the notation from Construction~$\ref{cons-flagregular}$ and from Lemma~\ref{lem:regmapgraph}, and assuming that $\Gamma=\Cos(G,H,J)$ is neither $\K_2^{(\lambda)}$ nor $\C_n^{(\lambda)}$ for any $n\geq3$,
\begin{enumerate}
\item[(a)] $\Cos(G,H,J,W)$ gives rise to a $G$-flag-regular map $\FRM(G,x,y,z)$ with valency $k$, edge-multiplicity $\lambda$ and face-length $m$, such that $(\alpha,e, f)$ is an incident triple, where $f=W\in F$,

\item[(b)] and, for each $g\in G$, the image of $C(W)$ under $g$ is the boundary cycle of the face $f^g$, and $G_{f^g}=G_{C(W)g}=W^g$; moreover $C(W)$ is a regular $m$-cycle and either $[C(W)] =\C_{m/\lambda'}^{(\lambda')}$ with $m/\lambda'\geq3$, or  $[C(W)] =\K_{2}^{(m)}$ with $m=2\lambda'\geq4$.

\item[(c)] The kernels of the $G$-actions on $V$ and $V\cup F$ are $G_{(V)}=\l a\r\cap \l a^z\r=\l a^k\r\cong\ZZ_\lambda$ and $G_{(V\cup F)}=\l a\r\cap \l b\r\cong \ZZ_{\lambda'}$. In particular, $\FRM(G,a,z)$ is a circular embedding of $\Gamma$ if and only if $G$ acts faithfully on $V\cup F$, or equivalently, if and only if $\lambda'=1$;

\item[(d)]
$x$ interchanges the two faces $f,f^x$ incident with the edge $e=J$, while $z$ fixes each of $f$ and $f^x$.

\end{enumerate}
\end{proposition}

\proof
By Lemma~\ref{lem:regmapgraph}, the underlying graph $\Gamma=(V,E,\bfI)$ is $\Cos(G,H,J)$ with valency $k$, edge-multiplicity $\lambda$, and by assumption  $\Gamma$ is neither $\K_2^{(\lambda)}$ nor $\C_n^{(\lambda)}$ for any $n$. To see that $\Cos(G,H,J,W)$ gives rise to a map with face set $F=[G:W]$, we first identify for each right coset of $W$ the set of edges incident with it. Since $G$ is transitive on $F$ it is sufficient to consider  $f=W$, and we note that the stabiliser $G_f=W$ in this action.
Since $H,J, W$ have pairwise non-empty intersections, $(\alpha,e,f)$ is an incident triple.
We claim that the set $E(f)$ of edges incident with $f$ in $\Cos(G,H,J,W)$ is precisely the set of edges in the cycle $C(W)$ of Lemma~\ref{lem:regmapgraph}(d). Suppose that the edge $Ju\in E(f)$, that is to say, $Ju\cap W\ne \emptyset$. Then $W$ contains at least one of the elements of
$Ju=\{u,xu,zu,zxu\}$. Since $z\in W$ this is equivalent to the condition that $W$ contains one of $u$ or $xu$. Also, since $x\in J$ we have $Ju=Jxu$, and replacing $u$ by $xu$ if necessary we may assume that $W$ contains $u$. Then $u=b^i$ or $zb^i$ for some $i$, and since $Jzb^i=Jb^i$ (as $z\in J$), it follows that $Ju=Jb^i$ for some $i$, proving the claim.

Now $x\in J\cap H=G_{e, \alpha}$, and by Lemma~\ref{lem:regmapgraph}(b), $x\not\in W=G_f$, and hence the face $f^x$ is incident with $e$ and distinct from $f$.
Each face is of the form $f^g=Wg$ for some $g\in G$. Suppose that $f^g$ is incident with $e$. Then by definition, $Wg\cap J\ne\emptyset$. If $Wg$ contains $1$ of $z$ then $Wg=W$ (since $z\in W$), and hence $f^g=f$. Otherwise $Wg$ contains $x$ or $zx$, and since $Wzx=Wx$ (as $z\in W$), this implies that $Wg=Wx$ and hence $f^g=f^x$. Thus $e$ is incident with precisely two faces, and this is true for all edges since $G$ is edge-transitive. The element $z\in W\cap J=G_{f,e}$ and hence $z$ fixes each of the two faces $f, f^x$ incident with $e$. This proves part (d).

Now the index 2 subgroup $\l a\r$ of $H=G_\alpha$ permutes cyclically and transitively the sets $E(\alpha)$ and $F(\alpha)$ of edges and faces incident with $\alpha$, respectively. Since $f$ is one of these faces and $f$ is incident with $e$, we may assume that $f^a$ is the second face $f^x$ incident with the arc $(\alpha, e)$.
Let $e_-:=e^{a^{-1}}$. Then $e_- \in E(\alpha)$ and $(e_-,e)$ is a $2$-arc on the boundary cycle $C(f)$ of $f$.
To identify the edge-sequence $C(f)$, note that $z\in W\cap J = G_{f,e}$ and $z$ interchanges the two vertices $\alpha$ and $\beta$ incident with $e$. Hence the action of $z$ on $C(f)$ is a reflection in the edge $[\alpha,e,\beta]$. Similarly $y\in W\cap H=G_{f,\alpha}$ and so the action of $y$ on $C(f)$ is a reflection in the vertex $\alpha$. It follows that the product $b=zy\in G_f$ acts on $C(f)$ as a rotation, sending $[\beta,e,\alpha]$ to $[\alpha,e_-, \beta^y]$. It follows that the edge $e_-$ is equal to $e^b=Jb$, and that $C(f)$ is the sequence $C(W)$ of Lemma~\ref{lem:regmapgraph}. Thus, by Lemma~\ref{lem:regmapgraph}, $C(f)$ is an $m$-cycle with stabiliser equal to $W$, that is, $G_f=G_{C(f)}=W$, and the assertions of part (b) follow, and the face-length is $m$.

To complete the proof of part (a), we identify each cycle $C(f^g)$, for $g\in G$, with the boundary of a Euclidean disc $\D(f^g)$, and
then `sew' these discs together by identifying arcs on two such discs when they correspond to the same arc of the corresponding cycles. This defines a topological space  $\calS$ with the graph $\Ga$ embedded in it. Proof that $\calS$ is a surface is identical to that given in the proof of Proposition~\ref{p:rotary}. Thus we have an embedding $\FRM(G,x,y,z)$ of $\Gamma$ in $\calS$ which admits the group $G$, and since $G$ is arc-transitive and the stabiliser of the arc $(\alpha,e)$ is $H\cap J=\l x\r$ and $x$ interchanges the two faces incident with $e$, it follows that $G$ is transitive on the flags, and hence  $\FRM(G,x,y,z)$ is $G$-flag-regular.

It remains to prove part (c). By our assumptions on $\Gamma$ it follows from Lemma~\ref{lem:regmapgraph}(d) that $G_{(V)}=\l a\r\cap \l a^z\r=\l a^k\r\cong\ZZ_\lambda$. Since $\l a\r\leq G_\alpha$ acts transitively on $F(\alpha)$, we have
 $G_{(V\cup F)}\leq \l a^k\r\cap G_f\leq H\cap W = (\l a\r\cap\l b\r)\rtimes\l y\r$.
 Now $\l a\r\cap\l b\r\cong\ZZ_{\lambda'}$ by Construction~\ref{cons-flagregular}, and hence acts trivially on $V$ while $y$ is nontrivial on $V$. Hence $G_{(V\cup F)}\leq
 \l a\r\cap\l b\r$, and since $ \l a\r\cap\l b\r$ is contained in the cyclic normal subgroup $N$ of Lemma~\ref{lem:regmapgraph}(d), it follows that this group acts trivially on $F$ and hence equality holds and  $G_{(V\cup F)}=\l a\r\cap\l b\r\cong\ZZ_{\lambda'}$.  In particular, $\FRM(G,a,z)$ is a circular embedding of $\Gamma$ if and only if $G$ acts faithfully on $V\cup F$, or equivalently, if and only if $\lambda'=1$.
\qed

\section{Vertex-rotary maps}\label{sec:vrmaps}

In this section we make precise  the dichotomy for $G$-vertex-rotary maps discussed in Remark~\ref{def:rotary-maps}, arising from the different possible actions of an edge-stabiliser on the two faces incident with the edge. We show in Proposition~\ref{p:vrmaps1} that each $G$-vertex-rotary map arises from one of the two general constructions we gave in Section~\ref{sec:cons}. In particular Theorem~\ref{thm} follows immediately from Proposition~\ref{p:vrmaps1}. We draw attention to the discussion in  Remark~\ref{rem:birot} of the bi-rotary case.

\begin{proposition}\label{p:vrmaps1}
Assume that Hypothesis~$\ref{Hypo}$ holds for a locally finite $G$-vertex-rotary map $\calM=(V, E, F, \bfI)$ with $|V|\geq 3$ and $|F|\geq 3$. So $f, f'$ are the two faces incident with the edge $e=\l z\r$ and we have $f^a=f'$ (replacing $a$ by $a^{-1}$ if necessary) and
\[
z:\ (f,f')\to (f,f')\ \mbox{or}\ (f',f).
\]

\begin{enumerate}
\item[(a)] If $z$ interchanges $f$ and $f'$, then $m:=|az|$ is finite and $\calM$ is the $G$-rotary map $\RM(G,a,z)$ of Construction~$\ref{cons-rotary}$ with face-length $m$.

\item[(b)] If $z$ fixes each of $f$ and $f'$, then $\ell:=|zz^a|$ is finite and  $\calM$ is the $G$-bi-rotary map $\BRM(G,a,z)$  of Construction~$\ref{cons-bi-rotary}$ with face-length $2\ell$.
\end{enumerate}
\end{proposition}

\proof
By Definition~\ref{def:rotary}(b), the underlying graph $\Ga=(V, E, \bfI)$ is $G$-vertex-rotary, and by Lemma~\ref{lem:rotpair}, we may assume that $\Gamma=\Cos(G,H,J)$, where for the arc $(\a, e)$, $G_\a=H=\l a\r$ and $G_e=J=\l z\r$, and $(a,z)$ is a rotary pair for $G$. Also $G$ is regular on the arc set of $\Gamma$.  By Hypothesis~\ref{Hypo}, $|V|\geq3$, $|E|\geq3$, and $|a|=k\lambda\geq3$. This implies in particular that $H\ne H^z$, for if $H=H^z$ then $G=H\rtimes \l z\r$ and $|V|=|G:H|=2$ which is a contradiction.
 As discussed in Remark~\ref{def:rotary-maps}, $z$ leaves $\{f,f'\}$ invariant, and $f^a=f'$ (replacing $a$ by $a^{-1}$ if necessary).

Since $a$ fixes $\a$ but not $e$, there exists an edge $e_-\in E(\a)$ such that $e_-\ne e$ and $(e_-)^a=e$. Further, since $\l a\r$ permutes transitively the $k\lambda$ faces of $\calM$ incident with $\a$, it follows that $(e_-,e)$ is a path of length two on the boundary cycle $C(f)$ of $f$ and, for $e_+:=e^a$, $(e,e_+)=(e_-,e)^a$ is  a path of length two on the boundary cycle $C(f')$ of $f'$. If $e_-=[\g,e_-,\a]$, then $\g^a=\b$, where $e=[\a,e,\b]$.
We claim that $\g\ne\b$. If $\g=\b$, then as $G_\a=\l a\r$ is transitive on $E(\a)$, it follows that every edge incident with $\a$ is also incident with $\b$, and since $\Ga$ is connected, we conclude that $V=\{\a,\b\}$, which is a contradiction. Thus $\g\ne\b$, proving the claim. Since $\g\ne\b$ and $\g^a=\b$, it follows that $\b^a\ne\b$, and the edge $e_+=[\a,e_+,\b^a]=e^a\ne e$. If $e_-=e_+$ then $a^2$ fixes the edge $e_-$, and  since $\l a\r$ acts faithfully and transitively on the $k\lambda$ edges of $E(\a)$ this implies that $|a|=k\lambda=2$, which is a contradiction. Thus $e_-, e, e_+$ are three pairwise distinct edges.

(a)
Assume first that $z$ interchanges $f$ and $f'$. Then $f^{az}=(f^a)^z=(f')^z=f$, so $\l az\r\leq G_f$. Now $(e_-,e)$ is a path of length two on $C(f)$, and
$(e_-)^{az}=e^z=e$, while $e^{az}$ is $[\a, e, \b]^{az}=[\a, e_+,\b^a]^z = [\b,(e_+)^z,\b^{az}]$, and this edge is not equal to $e_-$ since $e_-$ is not incident to $\b$.
The only induced action on the cycle $C(f)$ with these properties is a cyclic permutation of the edges of $C(f)$. This implies that $\l az\r$ is transitive on the edges of $C(f)$; moreover we conclude that the edge sequence $C(f)=(e, e^{az},\dots, e^{(az)^{m-1}})$, where $m=|az|$, is the sequence $C(az)$ of Lemma~\ref{lem:cyc-vertrot}(a). Hence $C(f)$ is a regular $m$-cycle with stabiliser $G_{C(f)}=\l az\r$ and induced subgraph $\C_{m/\lambda'}^{(\lambda')}$ where $\lambda'=|\l a\r\cap\l az\r|$ divides $\gcd(m,\lambda)$ and $m/\lambda'\geq3$. Since $G_f$ leaves $C(f)$ invariant we have $G_f\leq \l az\r$,
and since  $az$ leaves $f$ invariant we have equality $G_f=G_{C(f)}=\l az\r$. We observe the the edges  in $C(f)$ are precisely those of the form $\l z\r(az)^i$, for some $i$, and these are precisely the cosets of $\l z \r$ having non-empty intersection with $\l az\r$. Similarly the vertices incident with an edge of $C(f)$ are precisely those of the form $\l a\r (az)^i$ or $\l a\r z(az)^i$ for some $i$, and again these are precisely the cosets of $\l a\r$  with non-empty intersection with $\l az\r$. It follows that $\calM$ is the coset configuration $\Cos(G,H,J,\l az\r)$ of Construction~\ref{cons-rotary} and that $\calM$ is $\RM(G,a,z)$.

(b) Assume now that $z$ fixes each of $f$ and $f'$. Since $f'=f^a$ it follows that $z^a$ fixes $f'$, and hence $G_{f'}$ contains $\l z,z^a\r\cong \D_{2\ell}$, where $\ell=|zz^a|$. Recall from above that $(e,e_+)$ is a path of length $2$ on $C(f')$, where $e_+$ is $[\a, e_+, \b^a]$ with $\b\ne \b^a = \l a\r za$. Now $z, z^a$ induce reflections of $C(f')$ in the distinct consecutive edges $e, e_+$ of $C(f')$, respectively. This means in particular that $z\ne z^a$, and that the product $zz^a$ induces a `two-step rotation' $e'_i\to e'_{i+2}$ of the cycle $C(f')=(e_0',e_1',\dots)$. The first edge $e$ is the edge
$[\l a\r z, e_0',\a \r ]$, and the second edge $e_+$ is the edge $[\l a\r, e_1',\l a\r z(zz^a)]$, of the cycle $C(zz^a)$ of Lemma~\ref{lem:cyc-vertrot}(b), noting that $\l a\r z(zz^a)=\l a\r za$.  From the action of $\l zz^a\r$ we see that $C(f')$ is the edge sequence $C(zz^a)=(e_0',e_1',\dots, e'_{2\ell-1})$  of Lemma~\ref{lem:cyc-vertrot}(b).
Thus   $C(f')$ is a regular $2\ell$-cycle with stabiliser $G_{C(f')}=\l z,z^a\r\cong \D_{2\ell}$ and induced subgraph $\C_{2\ell/\lambda''}^{(\lambda'')}$ where $\lambda''=|\l a\r\cap\l zz^a\r|$ divides $\gcd(\ell,\lambda)$ and $\ell/\lambda''\geq3$. Since $G_{f'}$ leaves $C(f')$ invariant we have $G_{f'}\leq \l z,z^a\r$,
and since  $z, z^a$ both leave $f'$ invariant we have equality $G_{f'}=G_{C(f')}=\l z,z^a\r$. It follows that the incidence configuration of vertices, edges and faces of $\calM$ is precisely that of $\Cos(G, H, J, \l z,z^a\r)$ in Construction~\ref{cons-bi-rotary}, and, by Proposition~\ref{p:bi-rotary}, that $\calM$ is the $G$-bi-rotary map
$\BRM(G,a,z)$.
\qed

\begin{remark}\label{rem:birot}
{\rm
In Proposition~\ref{p:vrmaps1}(b) we saw that $G_\a=\l a\r$ cyclically permutes the $k\ell$ edges of $E(\a)$, sending $e$ to $e_+$. Hence $G_\b=\l a^z\r$ cyclically permutes $E(\b)$, sending $e^z=e$ to $(e_+)^z$. Since, in Proposition~\ref{p:vrmaps1} (b), $z$ fixes the face $f'$, and since $(e_+)^z$  is incident to $\a^z=\b$, it follows that $(e_+)^z$ is the edge $e_{2\ell-1}'$ of $C(f')$.
The stabilisers $G_\a=\l a\r$ and $G_\b=\l a^z\r$ thus induce cyclic permutations of $E(\a)$ and $E(\b)$, and hence local orientations of the carrier surface, which disagree. As $G$ is edge-transitive, this occurs for the two vertices incident with any edge. Such maps arising from embeddings of simple graphs have been called \emph{bi-orientable}, and in the case where the full map automorphism groups are arc-regular the maps are called \emph{bi-rotary}, see \cite[p.28]{bi-rotary}. As we do not assume that $G$ is the full automorphism group of $\calM$, we shall refer  to the maps in Proposition~\ref{p:vrmaps1}(b) as \emph{$G$-bi-rotary.}
}
\end{remark}

\section{Explicit families of vertex-rotary graphs and vertex-rotary maps}\label{sec:examples}

In this section, we prove Corollary~\ref{rotary=bi-rotary} by constructing examples  of rotary maps and bi-rotary maps with underlying graphs being hypercubes. We also give some examples arising from complete bipartite graphs.

\subsection{Hypercubes}

A simple $n$-dimensional hypercube $\bfQ_n$ is a Hamming graph $\bfH(n,2)$, which is a cartesian product of $n$ copies of $\K_2$, namely, $\K_2^{\square n}=\bfQ_n=\bfH(n,2)$.
We are interested in $\lambda$-extenders $\bfQ_n^{(\lambda)}$ with $\lambda\geqslant1$, and we construct such graphs as coset graphs.

\begin{definition}\label{def:cube}
{\rm
Define a group $A=\ZZ_2^n{:} \D_{2n\lambda}$ as follows:
\[\begin{array}{l}
A=\l v_0,v_1,\dots,v_{n-1}\r\rtimes(\l a\r{:}\l x\r),
\end{array}\]
where $|v_i|=2=|x|$ and $|a|=n\lambda$ such that
\[
\mbox{$a^x=a^{-1}$, $v_i^a=v_{i+1}$, and $v_i^x=v_{n-i}$},
\]
reading the subscripts modulo $n$.
In particular, $v_0x=xv_0$.
Let $z=v_0$, and let
\[
\begin{array}{ll}
X=\l v_0,v_1,\dots,v_{n-1}\r\rtimes \l a\r &=\ZZ_2^n{:}\ZZ_{n\lambda},\\
Y=(\l v_0v_1,v_1v_2,\dots, v_{n-2}v_{n-1}\r\rtimes\l a\r)\rtimes\l zx\r &=(\ZZ_2^{n-1}{:}\ZZ_{n\lambda}){:}\ZZ_2.
\end{array}
\]
}
\end{definition}

Then $X=\l a,z\r$, and $Y=\l a,zx\r$, so $(a,z)$ is a rotary pair for $X$, and $(a,zx)$ is a rotary pair for $Y$.

\begin{lemma}\label{coset-cube}
Let $y=ax$, and let $H=\l x,y\r=\l a\r{:}\l x\r=\D_{2n\lambda}$, and $J=\l x,z\r=\ZZ_2^2$. Then $(x,y,z)$ is a flag-regular triple for $A=\l x,y,z\r$, and $\bfQ_n^{(\lambda)}\cong \Cos(A,H,J)\cong\Cos(X,\l a\r,\l z\r)\cong \Cos(Y,\l a\r,\l zx\r)$, with  valency $n$ and edge-multiplicity $\lambda$.
\end{lemma}

\proof
It is easy to see that $\l x,y,z\r=A$.
Also, by Definition~\ref{coset-cube}, $xz=zx$, $z\not\in H=\l x,y\r$, and $A$ is finite, and hence $(x,y,z)$ is a fag-regular triple for $A$. Further, all the conditions of Construction~\ref{def-m-graph} hold for the subgroups $H, J$ of the group $A$.
Thus $\Ga=\Cos(A,H,J)$ is well defined.
By Theorem~\ref{thm:cosetgraph}, $\Ga$ is connected and $A$-arc-transitive, with vertex set $V=[A:H]$ and edge set $E=[A:J]$. Let $\a$ and $e$ be the vertex and the edge corresponding to $H$ and $J$, respectively.
The normal subgroup $N:=\l v_0,v_1,\dots,v_{n-1}\r$ of $A$ acts regularly on $V$ and hence is a complete set of coset representatives for $V=[A:H]$. Since $z\in J=G_e$, the vertex $\b=\a^z=Hz$ is incident with $e$ and lies in the $H$-orbit $\Gamma(\a)=\beta^H = \{ Hv_i\mid 0\leq i\leq n-1\}$ of size $n$.
Further $G_{\a,\b}=H\cap H^z=\l a^n\r\cong\ZZ_\lambda$,  and hence $\Gamma$ has valency $n$ and edge-multiplicity $\lambda$.

Now $A_\a=H=\l a\r{:}\l x\r$, and $X_\a=\l a\r=Y_\a$.
Note that both $\l v_0,v_1,\dots,v_{n-1}\r$ and $\l v_0v_1,v_1v_2,\dots, v_{n-2}v_{n-1}\r{:}\l zx\r$ are transitive on the vertex set $[A:H]$.
It follows that $X$ and $Y$ are arc-regular on $\Ga$.
Therefore, by Theorem~\ref{thm:cosetgraph2}, $\Ga$ is isomorphic to the coset graph $\Cos(X,\l a\r,\l z\r)$ of $X$, and the coset graph $\Cos(Y,\l a\r,\l zx\r)$ of $Y$,  so $\Ga=\Cos(A,H,J)\cong\Cos(X,\l a\r,\l z\r)\cong \Cos(Y,\l a\r,\l zx\r)$.

The group $L$ of Construction~\ref{def-m-graph} is $L=\l a^n,z\r = \l a^n\r \times\l z\r$ of order $2\lambda$. Hence by Theorem~\ref{thm:cosetgraph}(d), $\Gamma$ is an $(A,\lambda)$-extender of the base graph $\Gamma'=\Cos(A,H,L)\cong\SimpCos(A, H, HzH)$.  The simple base graph $\SimpCos(A, H, HzH)$ is easily recognised as the Cayley graph of $N\cong \ZZ_2^n$ with the identity $1_N$ adjacent to $\{v_0,\dots, v_{n-1}\}$; this is the cube of dimension $n$, namely, $\bfQ_n\cong\K_2^{\square n}$.
\qed

Using the notation defined in Definition~\ref{def:cube} and Lemma~\ref{coset-cube},
$(x,y,z)$ is a flag-regular triple for $A$ which, by Proposition~\ref{p:flagregular},  gives rise to a flag-regular map as defined in Construction~\ref{cons-flagregular}:
\[
\mbox{$\calM=\FRM(A,x,y,z)$, where $y=ax$.}
\]
The underlying graph of this map is $\Gamma:=\Cos(A,H,J)=\bfQ_n^{(\lambda)}$ by Lemma~\ref{coset-cube}.
Furthermore,
\begin{itemize}
\item[(a)] $(a,z)$ is a rotary pair for the group $X$, and
\item[(b)] $(a,zx)$ is a rotary pair for the group $Y$.
\end{itemize}
So we may construct vertex-rotary maps  by Constructions~\ref{cons-rotary} and~\ref{cons-bi-rotary} in each of these cases. We show that $\calM$ yields a map of each type when considered for these proper subgroups.

\begin{lemma}\label{rotary-bi-rotary}
With the notation above,
\[
\RegMap(A,x,y,z)=\BRM(X,a,z)=\RM(Y,a,zx)
\]
is a circular embedding of $\bfQ_n^{(\lambda)}$ with valency $n$, edge-multiplicity $\lambda$, and with simple face boundary cycles $\C_4$ of length $4$.
In particular the claims of Corollary~$\ref{rotary=bi-rotary}$ are valid.
\end{lemma}

\proof
By Proposition~\ref{p:flagregular}, $\calM=\RegMap(A,x,y,z)$ is an $A$-flag-regular map, and is defined as the coset configuration $\Cos(A, H, J, W)$ where $W=\l y,z\r$. Now $(\alpha, e, f)$ is an incident triple where $\alpha=H, e=J$ and $f=W$, and by Proposition~\ref{p:flagregular}(d), $xz$ interchanges the two faces $f, f^x$ incident with $e$ and $z$ fixes each of $f$ and $f^x$.

Now $X,Y<A$, and we showed in Lemma~\ref{coset-cube} that $X, Y$ are both arc-regular on the underlying graph $\Gamma$ of $\calM$. Since $H\cap X= H\cap Y = \l a\r$, it follows that $\Gamma$ is $X$-vertex-rotary and $Y$-vertex-rotary, and we note that $J\cap X=\l z\r$ and $J\cap Y=\l zx\r$. Since $z$ fixes each of $f, f^x$ and $zx$ interchanges them, it follows from Theorem~\ref{thm} (see also Propositions~\ref{p:rotary} and~\ref{p:bi-rotary}) that $\calM$ is the $X$-bi-rotary map $\BRM(X,a,z)$ and the $Y$-rotary map $\RM(Y,a,zx)$.

By Construction~\ref{cons-bi-rotary} and Proposition~\ref{p:bi-rotary},  $\BRM(X,a,z)$ has face-length $2\ell$ where $\ell=|zz^a|$. Since $z=v_0$ and $z^a=v_1$ we have $\ell=2$. Further, $zz^a=v_0v_1$ is an involution and does not lie in $\l a\r$, so we have $\l a\r\cap\l zz^a\r=1$ and it follows from Proposition~\ref{p:bi-rotary}(c) that
$\BRM(G,a,z)$ is a circular embedding with face-length $4$.
In particular the claims of Corollary~$\ref{rotary=bi-rotary}$ are valid.
\qed

\medskip
We remark that the group $X$ (as an abstract group) with the rotary pair $(a,z)$ also defines a rotary embedding of $\bfQ_n^{(\lambda)}$ by Construction~\ref{cons-rotary}.

\begin{lemma}\label{cube-map-1}
With the notation defined above, $\RM(X,a,z)$ is an $X$-rotary embedding  of $\bfQ_n^{(\lambda)}$ with valency $n$, edge-multiplicity $\lambda$, and face-length $2n\lambda/\gcd(2,\lambda)$.
Further, the face boundary cycles of $\calM$ are regular cycles with induced subgraphs
$\C_{2n}^{({\lambda/\gcd(2,\lambda)})}$.
\end{lemma}

\proof
By Proposition~\ref{p:rotary}\,(b), $\calM$ has face-length $m:=|az|$.
Modulo the normal subgroup $N:=\l v_0,v_1,\dots,v_{n-1}\r$, we see that $|Naz|=n$, and $(az)^n=va^n$ where $v=v_0v_1\dots v_{n-1}$. Since $|a^n|=\lambda$, $|v|=2$, and $a^n$ commutes with $v$, it follows that $|(az)^n| = 2\lambda/\gcd(2,\lambda)$, and hence that  $m=2n\lambda/\gcd(2,\lambda)$. Note that $2n$ divides $m$, and it is straightforward to prove that $\l a\r\cap\l az\r =\l a^{2n}\r\cong \ZZ_{\lambda/\gcd(2,\lambda)}$. Thus by Proposition~\ref{p:rotary}, the boundary cycle $C(f)$ is regular with induced subgraph as claimed.
\qed

\subsection{Complete bipartite graphs}

We give one more explicit family of examples of rotary embeddings.

Let $\Gamma=\K_{n,n}^{(\lambda)}$, where $n$ is odd, $\lambda$ is even, $\lambda > 2$,  and $\gcd(\lambda,n)=1$. Write $\lambda = 2\mu$ with $\mu>1$, and note that $\gcd(\lambda,\mu+1)=1$ if $\mu$ is even, and $\gcd(\lambda,\mu+2)=1$ if $\mu$ is odd.
Define a group:
\[
G=(\l b\r\times\l c_1\r\times \l c_2\r){:}\l z\r=(\ZZ_\lambda\times\ZZ_n^2){:}\ZZ_2,\ \mbox{where $|b|=\lambda$, $|c_1|=|c_2|=n$, and $|z|=2$,}
\]
such that
\[
(c_1,c_2)^z = (c_1,c_2^{-1}), \ \mbox{and\ $b^z=b^{\mu + \delta}$, where $\delta=1$ if $\mu$ is even, or $\delta=2$ if $\mu$ is odd.}
\]
Then $|G|=2\lambda n^2$.
Let
\[
a=bc_1c_2, \ \mbox{of order $|a|=\lambda n$.}
\]
We shall show that $(a,z)$ is a rotary pair for $G$ in Lemma~\ref{K(n,n):def}, and define a rotary embedding of $\K_{n,n}^{(\lambda)}$ in Lemma~\ref{K(n,n)-map:def}.


\begin{lemma}\label{K(n,n):def}
Using the notation defined above, $(a,z)$ is a rotary pair for the group $G$ and the $G$-vertex-rotary graph
$\Cos(G,\l a\r, \l z\r)\cong \K_{n,n}^{(\lambda)}$.
\end{lemma}

\proof
By Definition~\ref{def-rotarypair}, $(a,z)$ is a rotary pair for $G$ since $|z|=2$, $z\not\in\l a\r$, and
\[\begin{array}{rcll}
\l a,z\r&=&\l b,c_1c_2,z\r, &\mbox{since $\gcd(n,\lambda)=1$}\\
 		&=& \l b,c_1c_2,c_1c_2^{-1},z\r, &\mbox{as $(c_1c_2)^z=c_1c_2^{-1}$}\\
 		&=& \l b,c_1^2,c_1c_2^{-1},z\r, &\\
 &=& \l b,c_1^2,c_2^{-1},z\r, &\mbox{as $|c_1|=|c_2|=n$ is odd}\\
 &=& \l b,c_1,c_2,z\r=G.&
\end{array}\]

Let $\Ga=\Cos(G,\l a\r,\l z\r)=(V,E)$.
Then by Definition~\ref{def:rotary}(a) and Theorem~\ref{thm:cosetgraph}, $\Gamma$ is connected and $G$-vertex-rotary.
Also $|V|=|G|/|\l a\r|=2n$, $|E|=|G|/|\l z\r|=n^2\lambda$, and $G$ is regular on the arcs of $\Ga$.
{\color{black}
Since $\gcd(\lambda,n)=1$ and $n$ is odd, it follows that
\[\begin{array}{rcl}
\l a\r\cap\l a^z\r&=&\l bc_1c_2\r\cap\l b^{\mu+\delta}c_1c_2^{-1}\r\\
&=&(\l b\r\times\l c_1c_2\r)\cap(\l b^{\mu+\delta}\r\times\l c_1c_2^{-1}\r)\\
&=&(\l b\r\cap\l b^{\mu+\delta}\r)\times(\l c_1c_2\r\cap\l c_1c_2^{-1}\r)\\
&=&\l b^{\mu+\delta}\r\cong\ZZ_\lambda.
\end{array}\]
}

%

By Lemma~\ref{lem:rotpair2}, $\Gamma$ has valency $|\l a\r: \l a\r\cap\l a^z\r|=n$ and edge-multiplicity $|\l a\r\cap\l a^z\r|=\lambda$. The subgroups $K, L$ of Construction~\ref{def-m-graph} are $K=\l a\r\cap\l a^z\r=\l b^n\r$ and  $L=\l b^n, z\r$. Hence by Theorem~\ref{thm:cosetgraph} parts (c) and (d), $\Gamma=\Gamma_0^{(\lambda)}$, for the simple graph $\Gamma_0=\Cos(G, \l a\r, L)$, and by Theorem~\ref{thm:cosetgraph}(a) applied to $\Gamma_0$, $\Gamma_0$ also has valency $n$.  Now $G$ has an index $2$ (normal) subgroup $G^+:=\l b,c_1, c_2\r$ containing the vertex stabiliser $\l a\r$, and since $\Gamma$ is connected, it follows that $\Gamma$ and $\Gamma_0$ are bipartite. Thus
$\Gamma_0$ is a simple bipartite arc-transitive graph of valency $n$ and order $2n$, and hence  $\Ga_0\cong \K_{n,n}$ and $\Ga\cong \K_{n,n}^{(\lambda)}$.
\qed

The rotary pair $(a,z)$ defines a $G$-rotary map.

\begin{lemma}\label{K(n,n)-map:def}
Using the notation defined above, $\RM(G,a,z)$ (as defined in Construction~$\ref{cons-rotary}$) is a $G$-rotary embedding of $\K_{n,n}^{(\lambda)}$ with face-length $m:=|az|$, and the boundary cycle of each face is a regular $m$-cycle with induced subgraph $\C_{2n}^{\lambda'}$, where $m=2n\lambda'$ and $m, \lambda'$ are as in Table~$\ref{t:knn}$.
\end{lemma}

\begin{table}
\begin{center}
\begin{tabular}{lcc}
$\mu$ & $m$ & $\lambda'$   \\
\hline
$\mu\equiv 2\pmod{4}$& $\lambda n$ 	& $\lambda/2$\\
$\mu\equiv 0\pmod{4}$& $(\lambda/2) n$ 	& $\lambda/4$\\
$\mu\equiv \pm 1\pmod{6}$& $\lambda n$ 	& $\lambda/2$\\
$\mu\equiv 3\pmod{6}$& $(\lambda/3) n$ 	& $\lambda/6$\\
\hline
\end{tabular}
\caption{Table of values for $m, \lambda'$ for Lemma~\ref{K(n,n)-map:def} with $\lambda=2\mu$}\label{t:knn}
\end{center}
\end{table}

\proof
By Lemma~\ref{K(n,n):def}, $(a,z)$ is a rotary pair for $G$ and $\Ga=\Cos(G,\l a\r,\l z\r)\cong \K_{n,n}^{(\lambda)}$ is a $G$-vertex-rotary graph.
By Construction~\ref{cons-rotary} and Proposition~\ref{p:rotary},  $\calM=\RM(G,a,z)$ is a $G$-rotary map with face set $F=[G:\l az\r]$, face-length $m:=|az|$, and the boundary cycle of each face is a regular $m$-cycle with induced subgraph $\C_{m/\lambda'}^{\lambda'}$, where $\lambda'=|\l a\r\cap\l az\r|$. Thus we need to compute $m$ and $\lambda'$, and in particular to prove that $m/\lambda'=2n$. Let $\lambda=2\mu$ and $\delta$ be as above.

Note that $az=bc_1c_2z\notin \l b, c_1,c_2\r$, and hence $m=|az|$ is even. Also $\l az\r\cap\l a\r\leqslant\l a\r\leqslant\l b,c_1,c_2\r$, and
\[
aa^z=(az)^2=(bc_1c_2z)(bc_1c_2z)=bc_1c_2(zbc_1c_2z)=b^{\mu + \delta+1}c_1^2.
\]
Thus $m=|az|=2\cdot |aa^z|=2\cdot \lcm\{|b^{\mu+\delta+1}|, |c_1^2|\}$.
Further, $\l a\r\cap\l az\r=\l a\r\cap\l(az)^2\r$, and since $\gcd(\lambda,n)=1$, it follows that
\[
\l a\r\cap\l(az)^2\r=(\l b\r\times\l c_1c_2\r)\cap(\l b^{\mu+\delta+1}\r\times\l c_1^2\r)=
(\l b\r\cap \l b^{\mu+\delta+1}\r)\times(\l c_1c_2\r\cap\l c_1^2\r)=\l b^{\mu+\delta+1}\r.
\]
Thus $\lambda'=|\l a\r\cap\l az\r|=|b^{\mu+\delta+1}|$ is coprime to $n$, and so, since $n$ is even,
\[
m=2\cdot \lcm\{|b^{\mu+\delta+1}|, |c_1^2|\}=2\cdot\lcm\{\lambda',n\}=2\lambda'n.
\]
%
%
Finally we verify the entries of Table~\ref{t:knn}, noting that the value of $m$ follows from that of $\lambda'$ since $m=2\lambda'n$. Suppose first that $\mu$ is even. Then $\delta=1$ and $\gcd(\mu+\delta+1,\lambda)=\gcd(\mu+2, 2\mu)=\gcd(\mu+2,4)$. If $\mu\equiv 2\pmod{4}$, then we have $\lambda'=|b^2|=\lambda/2$, while if $\mu\equiv 0\pmod{4}$, then $\lambda'=|b^4|=\lambda/4$.  Suppose now that $\mu$ is odd. Then $\delta=2$ and $\gcd(\mu+\delta+1,\lambda)=\gcd(\mu+3, 2\mu)=\gcd(\mu+3,6)$. If $\mu\equiv \pm 1\pmod{6}$, then $\lambda'=|b^2|=\lambda/2$, while if $\mu\equiv 3\pmod{6}$, then $\lambda'=|b^6|=\lambda/6$.
%
\qed

\end{document}